\documentclass[11pt,a4paper]{article}%
\usepackage{cite}
\usepackage[utf8]{inputenc}
\usepackage[T1]{fontenc}
\usepackage[top=1cm, bottom=3cm, left=2.5cm, right=2.5cm, headheight=0.5cm,
footskip=60pt]{geometry}
\usepackage{amsmath}
\usepackage{indentfirst}
\usepackage{amscd}
\usepackage{amsthm}
\usepackage{amssymb}
\usepackage{graphicx}
\usepackage{latexsym}
\usepackage{color}
\usepackage{amsfonts}%
\providecommand{\U}[1]{\protect\rule{.1in}{.1in}}
\newtheorem{theorem}{Theorem}

\newtheorem{definition}[theorem]{Definition}

\newtheorem{lemma}[theorem]{Lemma}

\newtheorem{proposition}[theorem]{Proposition}
\newtheorem{remark}[theorem]{Remark}

\numberwithin{equation}{section}
\numberwithin{theorem}{section}
\begin{document}

\title{Global Existence of Weak Solutions \\for the Anisotropic Compressible Stokes System}
\author{D. Bresch\thanks{Univ. Grenoble Alpes, Univ. Savoie Mont-Blanc, CNRS, LAMA,
Chamb\'ery, France; didier.bresch@univ-smb.fr}, \hskip1cm C. Burtea
\thanks{Universit\'e Paris Diderot UFR Math\'ematiques Batiment Sophie
Germain, Bureau 727 8 place Aur\'elie Nemours, 75013 Paris;
cburtea@math.univ-paris-diderot.fr}\,}
\maketitle

\centerline{Dedicated to the memory of Genevi\`eve Raugel}

\begin{center}
{\small \textbf{Abstract} }
\end{center}

{\small In this paper, we study the problem of global existence of weak
solutions for the quasi-stationary compressible Stokes equations with an
anisotropic viscous tensor. This is done by comparing the limit of the
equations of the energies associated to a sequence of weak-solutions with the
energy equation associated to the system verified by the limit of the sequence
of weak-solutions. This allows us to construct a particular defect measure
associated to the pressure which yields compactness. By doing so we avoid the
use of the so-called effective flux. Using this new tool, we solve an open
problem namely global existence of solutions \`{a} la Leray for such a system
without assuming any restriction on the anisotropy amplitude. This provides a
flexible and natural method to treat compressible quasilinear Stokes systems
which are important for instance in biology, porous media, supra-conductivity
or other applications in the low Reynolds number regime. }



\section{Introduction}

\subsection{Presentation of the main result}

As explained in \cite{Li}, Chapter $8$, there are various motivations for the
study of quasi-stationary Stokes problem. On the one hand such a study may be
used to try to understand how to construct solutions of the compressible
Navier-Stokes system which exhibit persistent oscillations. On the other hand
this system naturally arrises either when dealing with flows in the low
Reynolds number regime, which it typically the case in porous media or biology
either as a mean field model for the motion of vortices in a superconductor in
the Ginzburg--Landau theory. There is a rather rich literature regarding the
mathematical study assuming isotropic diffusion: see for instance
\cite{BrJa1}, \cite{Ka}, \cite{Ka1}, \cite{Ka2}, \cite{Ma}, \cite{MaZh}, or
\cite{VaKa} for constant viscosity coefficients or \cite{BePi} for density
dependent viscosity coefficients. More complicated versions of the
quasi-stationary compressible Stokes system have been also analyzed in
\cite{BrMuZa}, \cite{BrNePe}, \cite{FrWe} and \cite{FrGoMa} in the multi-fluid setting.

Global existence of weak solutions for general anisotropic viscosities for
non-stationary compressible barotropic Navier-Stokes equations or even
quasi-stationary Stokes equations are open problems. Only recently a positive
result has been obtained by D. Bresch and P.--E. Jabin in \cite{BrJa} assuming
some restrictions on the shear and bulk viscosities. The result is not
straightforward to prove as the anisotropy introduces non-locality in the
compactness characterization process. This explains in some sense the new
method introduced by the authors in order to conclude compactness: propagation
of a non-local $L^{p}$-compactness module with appropriate time-evolving weights.

In this paper, we consider a very general form of the quasi-stationary
compressible Stokes equations:
\begin{equation}
\left\{
\begin{array}
[c]{l}%
\partial_{t}\rho+\operatorname{div}\left(  \rho u\right)  =0,\\
-\mathrm{div}\ {\tau}+a\nabla\rho^{\gamma}=f,
\end{array}
\right.  \label{ANISYS}%
\end{equation}
completed with an initial density distribution
\begin{equation}
\rho|_{t=0}=\rho_{0}\geq0. \label{densityini}%
\end{equation}
Above, $u$ stands for the fluid velocity field, $\rho$ is the fluid density
and $\tau$ represents the viscous stress tensor which is given by\footnote{We
use the convention of summation over repeated indices.}
\begin{equation}
\tau_{ij}(t,x,D(u))=A_{ijkl}(t,x)[D(u)]_{kl} \label{tau}%
\end{equation}
where $D(u)=(\nabla u+{}^{t}\nabla u)/2$ is the strain tensor and
\begin{equation}
A_{ijkl}=A_{ijkl}(t,x)\in W^{1,\infty}((0,T)\times{\mathbb{T}}^{3}) \label{A}%
\end{equation}
are given coefficients. Also, $a>0$ is a given constant. The classical
isotropic case is obtained by choosing%
\[
\left\{
\begin{array}
[c]{l}%
A_{iiii}=\left(  \mu+\lambda\right)  ,\\
A_{iijj}=\lambda\text{ for }i\not =j,\\
A_{ijij}=A_{ijji}=\frac{\mu}{2}\text{ for }i\not =j,\\
A_{ijkl}=0\text{ otherwise.}%
\end{array}
\right.
\]
The simplest case example of anisotropic viscous stress tensor is obtained
for
\begin{equation}
\left\{
\begin{array}
[c]{l}%
A_{1111}=\mu_{1},A_{2222}=\mu_{2},\text{ }A_{3333}=\mu_{3},\\
A_{ijkl}=0\text{ otherwise,}%
\end{array}
\right.  \label{Temam-Ziane}%
\end{equation}
case in which we have
\[
\operatorname{div}\tau=\partial_{11}u+\partial_{22}u+\mu\partial_{33}%
u\overset{not.}{=}\Delta_{\mu}u.
\]
The aim of this paper is to present a proof in the spirit of that of Lions for
the existence and the weak stability of solutions i.e. we introduce a
particular defect measure for the pressure which allows to control the
oscillation of an approximating sequence of solutions of system
\eqref{ANISYS}--\eqref{densityini}. Of course, the key point that allows to
account for anisotropy is that we are able to control this defect measure
without using the effective flux. For the reader's convenience we will present
a sketch of the proof in the next section in the case of the viscous tensor
given by $\left(  \text{\ref{Temam-Ziane}}\right)  $.

In order to obtain a satisfactory mathematical theory we need to further
assume the following hypothesis on the stress tensor $\tau$:
\begin{align}
\bullet &  A_{ijkl}=A_{ijlk}\text{ for all }i,j,k,l\text{ which allows us to
write that}\nonumber\\
&  \tau(t,x,D(u)):\nabla u=\frac{1}{2}\tau(t,x,D(u)):D(u)\label{H1}\\
\bullet &  D(u)\longmapsto\tau
(t,x,D(u)):D(u)\hbox{ to be weakly lower semi-continuous }\label{H2}\\
\bullet &  \hbox{There exists }c>0\hbox{ such that }\nonumber\\
&  \hskip2cmE=\int_{{\mathbb{T}}^{3}}\tau(t,x,D(u)):\nabla u\geq
c\int_{{\mathbb{T}}^{3}}|\nabla u|^{2}\label{H3}\\
\bullet &  \hbox{ The application }\mathcal{A}:v\mapsto-\mathrm{div}%
\,\tau(t,x,D(v))\nonumber\\
&  \hbox{ is a second  order invertible elliptic operator}\nonumber\\
&  \hbox{ such that }\mathcal{A}^{-1}\nabla\mathrm{div}\text{ is a bounded
operator from }L^{\frac{3}{2}-\delta}\left(  \mathbb{T}^{3}\right)  \text{
into }L^{\frac{3}{2}-\delta}\left(  \mathbb{T}^{3}\right)
\hbox{for some}\nonumber\\
&  \delta\in(0,1/2). \label{H4}%
\end{align}
\noindent We are now in the position of announcing our main result:

\begin{theorem}
\label{Main}Consider $f,\partial_{t}f\in L^{2}((0,T);L^{\frac{6}{5}}\left(
\mathbb{T}^{3}\right)  )$ and initial data $\rho_{0}$ satisfying
\[
\rho_{0}\geq0,\qquad0<M_{0}=\int_{\mathbb{T}^{3}}\rho_{0}<+\infty,\qquad
E_{0}=\int_{\mathbb{T}^{3}}\rho_{0}^{\gamma}\,dx<+\infty,\qquad\int
_{\mathbb{T}^{3}}f\left(  t\right)  dx=0,
\]
where $\gamma>1$ and assume that the viscous stress tensor $\tau$ given by
\eqref{tau} satisfies \eqref{H1}--\eqref{H4}. Then there exists a global weak
solution $(\rho,u)$ of the system \eqref{ANISYS} and \eqref{densityini} with
\[
\rho\in\mathcal{C}([0,T];L_{weak}^{\gamma}(\mathbb{T}^{3}))\cap L^{2\gamma
}((0,T)\times\mathbb{T}^{3}),\quad u\in L^{2}(0,T;H^{1}(\mathbb{T}%
^{3})\hbox{ with }\text{ }\int_{\mathbb{T}^{3}}u=0.
\]

\end{theorem}

A similar result can be obtained for the case of a bounded domain with
Dirichlet boundary condition: we have chosen periodic boundary conditions to
simplify the presentation. One of the most delicate points in proving Theorem
\ref{Main} is the stability of weak-solutions namely, given a sequence of
solutions $\left(  \rho^{\varepsilon},u^{\varepsilon}\right)  $ of $\left(
\text{\ref{ANISYS}}\right)  $ verifying uniformly the energy estimates and
therefore (at least on a subsequence) weakly converge to some $\left(
\rho,u\right)  $, show that $\left(  \rho,u\right)  $ is also a solution for
$\left(  \text{\ref{ANISYS}}\right)  $. Of course, the most difficult part is
to identify the pressure term in the limit i.e. to prove that $\lim\>
(\rho^{\varepsilon})^{\gamma}=\rho^{\gamma}$. Of course, the case $\gamma=1$
does not present this difficulty. This is the reason why we choose to focus
only on the "more nonlinear" cases $\gamma>1$.

\begin{remark}
As explained in \textrm{\cite{Li}}, including a force term in the momentum
equation which is of the form $\rho g$ say with $g\in L_{t,x}^{\infty}$ does
not always have a solution because one has the compatibility condition%
\[%
{\displaystyle\int_{\mathbb{T}^{3}}}
(\rho g+f)=0.
\]
Thus, if $g$ is a vectors with positive components this would imply that
$\rho=0$ for all times and this independently of the initial data.
\end{remark}

One limitation of our work seems to be the choice of the pressure function: we
cannot consider more general convex pressure laws other than $p\left(
\rho\right)  =a\rho^{\gamma}$, see Remark \ref{observation1}. Also, it seems
difficult to adapt the method presented in this paper to the non-stationary
Navier-Stokes system for a compressible fluid. Note that actually only one
result exists for this system in the case of anisotropic diffusion, see
\cite{BrJa}. Loosely speaking, the authors require that the "quantity of
anisotropy" that they allow in the system should be small compared to the
total viscosity $2\mu+\lambda$. Observe that we do not impose such restriction
for the quasi-stationary Stokes system. However we are able to treat a
stationary system that can be interpreted as an implicit discretization of the
full Navier-Stokes system, see Section \ref{Applications}.

\medskip

The rest of the paper is organized as follows:

\begin{itemize}
\item Section \ref{Newapproach} is dedicated to present the new defect measure
associated to the pressure and to show how it is possible to control it if
this is the case initially. Our result uses in a crucial manner compactness
properties on the velocity field in $L^{2}((0,T)\times{\mathbb{T}}^{3})$. For
the readers's convenience, we recall the classical approach due to
P.--L.~Lions and latter refined by E. Feireisl-A. Novotny-H.Petzeltova. In
particular, we explain why the anisotropic case seems to fall completely out
of such strategy (see also \cite{BrJa} for further discussions).
\end{itemize}

The rest of the paper is devoted to the proof of Theorem \ref{Main}. As it is
accustomed when dealing with the existence of weak solutions, the proof is
divided into two parts.

\begin{itemize}
\item In Section \ref{WeakStability} we define and investigate the stability
of a sequence of bounded-energy weak-solutions of the system $\left(
\text{\ref{ANISYS}}\right)  $. In Section \ref{Section_tools} we recall the
basic nonlinear analysis tools that allow us to render rigorous the formal
computations presented in Section \ref{Newapproach}. In Section
\ref{Aprioriestimates} we prove that bounded energy-weak-solutions enjoy
extra-integrability and time regularity properties, with respect to the basic
energy estimates, of course. More precisely it turns out that $\rho^{\gamma
}\in L_{t,x}^{2}$ and that $\partial_{t}u\in L^{1}\left(  0,T;L^{r}\left(
\mathbb{T}^{3}\right)  \right)  $ for some $r\in(1,3/2)$. In Section
\ref{Weakstability} we investigate the stability of a sequence of bounded
energy weak-solutions $\left(  \rho^{\varepsilon},u^{\varepsilon}\right)  $
satisfying uniformly the energy estimates. It turns out that comparing the
limit of the energy associated to each solution $\left(  \rho^{\varepsilon
},u^{\varepsilon}\right)  $ with the energy of the system verified by $\left(
\rho,u\right)  =\lim\left(  \rho^{\varepsilon},u^{\varepsilon}\right)  $ we
obtain an identity that involves a defect measure associated to the pressure.
The stability result, interesting in itself is formalized in Theorem
\ref{main}, and it can be adapted to construct solutions for the system
$\left(  \text{\ref{ANISYS}}\right)  $.

\item In section \ref{Construction} we construct weak-solutions for the system
$\left(  \text{\ref{ANISYS}}\right)  $. More precisely, we propose an
approximate model that depends on two parameters such that, at least formally,
system $\left(  \text{\ref{ANISYS}}\right)  $ is obtained by a limit process
by making the parameters tend to zero. We show that we can construct solutions
by a classical fixed-point argument for the approximate system. Moreover, we
show that the solutions verify uniform bounds with respect to the parameters
introduced such that we are able to pass to the limit in a sequence of
solutions and show that the limiting object is a solution of for the system
$\left(  \text{\ref{ANISYS}}\right)  $ and thus achieving the proof of Theorem
\ref{Main}.

\item Finally, in Section \ref{Applications} \ we discuss some extents of our
method of proof to other systems.
\end{itemize}

\subsection{Formal approach to control the defect measure associated to the
pressure in a simplified case\label{Newapproach}}

To be understandable for the reader, let us present formally on a simple
example why the classical approach to control defect measures fails to apply
in the case of anisotropic viscosities and how our new way to proceed provides
a flexible method for Stokes type systems. More precisely, let us consider
$\left(  \rho^{\varepsilon},u^{\varepsilon}\right)  $ a sequence of solutions
for the following system.
\begin{equation}
\left\{
\begin{array}
[c]{l}%
\partial_{t}\rho^{\varepsilon}+\operatorname{div}\left(  \rho^{\varepsilon
}u^{\varepsilon}\right)  =0,\\
-\Delta_{\mu}u^{\varepsilon}+\nabla((\rho^{\varepsilon})^{\gamma})=f
\end{array}
\right.  \label{equations}%
\end{equation}
where
\[
\Delta_{\mu}=\mu_{1}\partial_{11}+\mu_{2}\partial_{22}+\mu_{3}\partial_{33}%
\]
with $\mu_{1},\mu_{2},\mu_{3}>0$ which may be different. Assume
\[
\Vert u^{\varepsilon}\Vert_{L^{2}(0,T;H^{1}(\mathbb{T}^{3}))}+\Vert
\rho^{\varepsilon}\Vert_{L^{2\gamma}((0,T)\times\mathbb{T}^{3})}+\Vert
\rho^{\varepsilon}\Vert_{L^{\infty}(0,T;L^{\gamma}(\mathbb{T}^{3}))}\leq
C<+\infty
\]
where $C$ does not depend on $\varepsilon$ weak solutions of \eqref{equations}
and assume that
\[
\{u^{\varepsilon}\}_{\varepsilon}\text{ is compact in }L^{2}((0,T)\times
{\mathbb{T}}^{3}).
\]
We denote $\left(  \rho,u\right)  $ the weak limit and, using classical
functional analysis arguments it is not hard to see that we have%
\begin{equation}
\left\{
\begin{array}
[c]{l}%
\partial_{t}\rho+\operatorname{div}\left(  \rho u\right)  =0,\\
-\Delta_{\mu}u+\nabla(\overline{\rho^{\gamma}})=f.
\end{array}
\right.  \label{equationslim}%
\end{equation}
for some function $\overline{\rho^{\gamma}}\in L^{2}((0,T)\times{\mathbb{T}%
}^{3})$. Of course, the main difficulty is to prove that $\overline
{\rho^{\gamma}}=\rho^{\gamma}$ and therefore to be able to characterize the
possible defect measures.

\begin{remark}
Throughout the paper we denote the weak limit of a sequence $\left(
a^{\varepsilon}\right)  _{\varepsilon>0}$ by $\bar{a}$.
\end{remark}

\noindent\textit{Classical approach to control defect measures.} As mentioned
in \cite{BrJa}, the usual method for isotropic viscosities (namely $\mu
_{1}=\mu_{2}=\mu_{3}=\mu$) is based on the careful analysis of the defect
measures
\[
\mathrm{dft}[\rho^{\varepsilon}-\rho](t)=\int_{\mathbb{T}^{3}}(\overline
{\rho\log\rho})(t)-\rho\log\rho(t))\,dx.
\]
More precisely, we can write the two equations
\begin{equation}
\partial_{t}(\rho\log\rho)+\mathrm{div}(\rho\log\rho u)+\rho\mathrm{div}u=0
\label{equ1}%
\end{equation}
and
\begin{equation}
\partial_{t}(\overline{\rho\log\rho})+\mathrm{div}(\overline{\rho\log\rho
}u)+\overline{\rho\mathrm{div}u}=0 \label{equ2}%
\end{equation}
Note that if $\rho\in L^{2}((0,T)\times{\mathbb{T}}^{3})$ then using the
uniform bound on $u\in L^{2}(0,T;H^{1}({\mathbb{T}}^{3}))$, we have
$\rho\,\mathrm{div}u\in L^{1}((0,T)\times{\mathbb{T}}^{3})$ and therefore the
third quantity is well defined. At this level comes the so called effective
flux comes into play. More precisely, Lions \cite{Li2} in $^{\prime}93$ (see
also D. Serre \cite{Serre} for the $1d$ case) observes that the following
quantity
\[
F^{\varepsilon}=p(\rho^{\varepsilon})-\mu\mathrm{div}u^{\varepsilon}%
\]
enjoys the following compactness property:
\begin{equation}
\lim_{\varepsilon\rightarrow0}\int_{0}^{T}\int_{\mathbb{T}^{3}}(p(\rho
^{\varepsilon})-\mu\mathrm{div}u^{\varepsilon})b(\rho^{\varepsilon}%
)\varphi=\int_{0}^{T}\int_{\mathbb{T}^{3}}(\overline{p(\rho)}-\mu
\mathrm{div}u)\overline{b(\rho)}\varphi. \label{flux}%
\end{equation}
This is important at it provides a way to express $\overline{\rho
\mathrm{div}u}$ in terms of $\rho\mathrm{div}u$ and an extra term which is
signed. Substracting the two equations $\left(  \text{\ref{equ1}}\right)  $
and $\left(  \text{\ref{equ2}}\right)  $ and using the important property of
the effective flux $\left(  \text{\ref{flux}}\right)  $, one gets that
\[
\partial_{t}(\overline{\rho\log\rho}-\rho\log\rho)+\mathrm{div}((\overline
{\rho\log\rho}-\rho\log\rho)u)=\frac{1}{\mu}(\overline{p(\rho)}\rho
-\overline{p(\rho)\rho})
\]
and using the monotonicity of the pressure, one may deduce that
\[
\mathrm{dft}[\rho^{\varepsilon}-\rho](t)\leq\mathrm{dft}[\rho^{\varepsilon
}-\rho](0).
\]
On the other hand, the strict convexity of the function $s\mapsto s\log s$
with $s\geq0$ implies that $\mathrm{dft}[\rho^{\varepsilon}-\rho](t)\geq0$. If
initially this quantity vanishes, it then vanishes at every time. The
commutation of the weak convergence with a strictly convex function yields
compactness of $\{\rho^{\varepsilon}\}_{\varepsilon}$ in $L^{1}((0,T)\times
{\mathbb{T}}^{3})$.

Assuming anisotropic viscosities $\mu_{1}=\mu_{2}\not =\mu_{3}$, the effective
flux property reads
\[
\overline{\rho\mathrm{div}u}-\rho\mathrm{div}u=\frac{1}{\mu_{1}}[\rho
\overline{A_{\nu}\rho^{\gamma}}-\overline{\rho A_{\nu}\rho^{\gamma}}]
\]
with some non-local anisotropic operator $A_{\nu}=(\Delta-(\mu_{3}-\mu
_{1})\partial_{z}^{2})^{-1}\partial_{z}^{2}$ where $\Delta$ is the total
Laplacian in terms of $(X,z)$ with variables $X=(x,y)$ and $z$. Unfortunately,
we are loosing the sign of the right-hand side. This explains why the
anisotropic case seems to fall completely out the theory developed by P.--L.
Lions \cite{Li} and E. Feireisl, A. Novotny and H. Petzeltova \cite{FeNoPe}.
The first positive answer has been given by D. Bresch and P.-E. Jabin in
\cite{BrJa} for the compressible Navier-Stokes equations developing an other
way to characterize compactness in space on the density: it involves a
non-local compactness criterion with the introduction of appropriate weights.
It allows them to obtain a positive answer assuming the viscosity coefficient
$\mu_{1},\mu_{2},$ $\mu_{3}$ to be close enough.

\medskip

\noindent\textit{New approach to control defect measures in the Stokes
regime.} Our new approach is based on the careful analysis of the defect
measures
\[
\mathrm{dft}[\rho^{\varepsilon}-\rho](t)=\int_{\mathbb{T}^{3}}\Bigl((\overline
{\rho^{\gamma}})(t)-\rho^{\gamma}(t)\Bigr)^{1/\gamma}\,dx.
\]
The main idea here is to write the equation related to the energy which will
not use the effective flux expression but is related to the viscous
dissipation in the Stokes regime. More precisely, let us observe that the
pressure \ verifies the following equation :%
\[
\partial_{t}\left(  \rho^{\varepsilon}\right)  ^{\gamma}+\operatorname{div}%
\left(  \left(  \rho^{\varepsilon}\right)  ^{\gamma}u^{\varepsilon}\right)
+\left(  \gamma-1\right)  \left(  \rho^{\varepsilon}\right)  ^{\gamma
}\operatorname{div}u^{\varepsilon}=0
\]
which rewrites%
\[
\partial_{t}\left(  \rho^{\varepsilon}\right)  ^{\gamma}+\gamma
\operatorname{div}\left(  \left(  \rho^{\varepsilon}\right)  ^{\gamma
}u\right)  -\left(  \gamma-1\right)  u^{\varepsilon}\nabla\left(
\rho^{\varepsilon}\right)  ^{\gamma}=0.
\]
We observe that with the aid of the second equation of $\left(
\text{\ref{equations}}\right)  $ we may write that%
\[
\partial_{t}\left(  \rho^{\varepsilon}\right)  ^{\gamma}+\gamma
\operatorname{div}\left(  \left(  \rho^{\varepsilon}\right)  ^{\gamma
}u\right)  -\left(  \gamma-1\right)  u^{\varepsilon}\Delta_{\mu}%
u^{\varepsilon}=(\gamma-1)u^{\varepsilon}f
\]
which can be put under the following form%
\begin{align}
&  \partial_{t}\left(  \rho^{\varepsilon}\right)  ^{\gamma}+\gamma
\operatorname{div}\left(  \left(  \rho^{\varepsilon}\right)  ^{\gamma
}u\right)  -\left(  \gamma-1\right)  \Delta_{\mu}\left(  \frac{\left\vert
u^{\varepsilon}\right\vert ^{2}}{2}\right) \nonumber\\
&  \hskip2cm=-\left(  \gamma-1\right)  \nabla_{\mu}u^{\varepsilon}:\nabla
_{\mu}u^{\varepsilon}+(\gamma-1)u^{\varepsilon}f, \label{equation2}%
\end{align}
where we use the notation%
\[
\nabla_{\mu}=\left(  \mu_{1}^{\frac{1}{2}}\partial_{1},\mu_{2}^{\frac{1}{2}%
}\partial_{2},\mu_{3}^{\frac{1}{2}}\partial_{3}\right)  .
\]
Of course, we used that
\[
\partial_{jj}u_{i}u_{i}=\partial_{jj}\left(  \frac{\left(  u_{i}\right)  ^{2}%
}{2}\right)  -(\partial_{j}u_{i})^{2}%
\]
Assuming that
\[
\left(  u^{\varepsilon}\right)  _{\varepsilon>0}\text{ is compact in }%
L^{2}((0,T)\times{\mathbb{T}}^{3})
\]
by passing to the limit in $\left(  \text{\ref{equation2}}\right)  $ we obtain
that%
\begin{align}
&  \partial_{t}\overline{\rho^{\gamma}}+\gamma\operatorname{div}\left(
\overline{\rho^{\gamma}}u\right)  -\left(  \gamma-1\right)  \Delta_{\mu
}\left(  \frac{\left\vert u\right\vert ^{2}}{2}\right) \nonumber\\
&  \hskip2cm=-\left(  \gamma-1\right)  \overline{\nabla_{\mu}u:\nabla_{\mu}%
u}+(\gamma-1)fu^{\varepsilon}. \label{Imp1}%
\end{align}
In the following we will apply the same recipe to the limiting function
$\left(  \rho,u\right)  $. Indeed, from $\left(  \text{\ref{equationslim}%
}\right)  $ one can deduce that
\begin{align*}
\partial_{t}\rho^{\gamma}+\gamma\operatorname{div}\left(  \rho^{\gamma
}u\right)   &  =\left(  \gamma-1\right)  u\cdot\nabla\rho^{\gamma}\\
&  =\left(  \gamma-1\right)  u\cdot\nabla(\rho^{\gamma}-\overline{\rho
^{\gamma}})-\left(  \gamma-1\right)  u\cdot\nabla\overline{\rho^{\gamma}}\\
&  =-\left(  \gamma-1\right)  u\cdot\nabla(\overline{\rho^{\gamma}}%
-\rho^{\gamma})-\left(  \gamma-1\right)  u\cdot(\Delta_{\mu}u+f)
\end{align*}
which rewrites%
\begin{align}
&  \partial_{t}\rho^{\gamma}+\gamma\operatorname{div}\left(  \rho^{\gamma
}u\right)  +\left(  \gamma-1\right)  u\nabla(\overline{\rho^{\gamma}}%
-\rho^{\gamma})-\left(  \gamma-1\right)  \Delta_{\mu}\left(  \frac{\left\vert
u\right\vert ^{2}}{2}\right) \nonumber\\
&  =-\left(  \gamma-1\right)  \nabla_{\mu}u:\nabla_{\mu}u+(\gamma-1)fu.
\label{Imp2}%
\end{align}
Let us consider the difference between $\left(  \text{\ref{Imp1}}\right)  $
and $\left(  \text{\ref{Imp2}}\right)  $ in order to write that%
\begin{align*}
&  \partial_{t}\left(  \overline{\rho^{\gamma}}-\rho^{\gamma}\right)
+\gamma\operatorname{div}\left(  \left(  \overline{\rho^{\gamma}}-\rho
^{\gamma}\right)  u\right)  -\left(  \gamma-1\right)  u\nabla\left(
\overline{\rho^{\gamma}}-\rho^{\gamma}\right) \\
&  =-\left(  \gamma-1\right)  \left(  \overline{\nabla_{\mu}u:\nabla_{\mu}%
u}-\nabla_{\mu}u:\nabla_{\mu}u\right)  .
\end{align*}
which we put under the form%
\begin{align}
&  \partial_{t}\left(  \overline{\rho^{\gamma}}-\rho^{\gamma}\right)
+\operatorname{div}\left(  \left(  \overline{\rho^{\gamma}}-\rho^{\gamma
}\right)  u\right)  +\left(  \gamma-1\right)  \left(  \overline{\rho^{\gamma}%
}-\rho^{\gamma}\right)  \operatorname{div}u\label{Imp3}\\
&  =-\left(  \gamma-1\right)  \left(  \overline{\nabla_{\mu}u:\nabla_{\mu}%
u}-\nabla_{\mu}u:\nabla_{\mu}u\right)  .\nonumber
\end{align}
At this point we observe that owing to the convexity of the pressure function,
we have that
\[
\overline{\rho^{\gamma}}\geq\rho^{\gamma}\text{ a.e.}%
\]
and
\begin{equation}
\overline{\nabla_{\mu}u:\nabla_{\mu}u}-\nabla_{\mu}u:\nabla_{\mu}u\geq0
\label{positive_mesure}%
\end{equation}
at least in the sense of measures. By multiplying $\left(  \text{\ref{Imp3}%
}\right)  $ with $\frac{1}{\gamma}\left(  \overline{\rho^{\gamma}}%
-\rho^{\gamma}\right)  ^{\frac{1}{\gamma}-1}$ we get that
\[
\partial_{t}\left(  \overline{\rho^{\gamma}}-\rho^{\gamma}\right)  ^{\frac
{1}{\gamma}}+\operatorname{div}\left(  \left(  \overline{\rho^{\gamma}}%
-\rho^{\gamma}\right)  ^{\frac{1}{\gamma}}u\right)  \leq0
\]
such that by integration and using $\left(  \text{\ref{positive_mesure}%
}\right)  $ we end up with
\[
\int_{0}^{T}\int\left(  \overline{\rho^{\gamma}}-\rho^{\gamma}\right)
^{\frac{1}{\gamma}}\leq T\int\left(  \overline{\rho^{\gamma}}-\rho^{\gamma
}\right)  _{|t=0}^{\frac{1}{\gamma}}.
\]
Therefore if we have compactness initially, we get compactness of the sequence
$\left(  \rho_{\varepsilon}\right)  _{\varepsilon\geq0}$. Of course all the
previous formal calculations have to be justified because of the weak
regularity and of possible vanishing quantity: this will be the subject of
Subsection \ref{weakstab}.

\begin{remark}
It is interesting to note that our new approach to get characterization of the
defect measure on the pressure sequence is related to the energy equation and
strongly uses the energy dissipation. We speculate that it has a physical
meaning in some sense.
\end{remark}

\begin{remark}
\label{observation1}Even though our method allows us to treat very general
anisotropies it does not seem to apply to general convex pressure laws
$p\left(  \rho\right)  $. If we let $H\left(  \rho\right)  $ be the potential
energy which is defined via%
\[
\rho H^{\prime}\left(  \rho\right)  -H\left(  \rho\right)  =p\left(
\rho\right)  ,
\]
then, we still have the identity%
\begin{align*}
&  \partial_{t}\left(  \overline{H\left(  \rho\right)  }-H\left(  \rho\right)
\right)  +\operatorname{div}\left(  \left(  \overline{H\left(  \rho\right)
}-H\left(  \rho\right)  \right)  u\right)  +\left(  \overline{p\left(
\rho\right)  }-p\left(  \rho\right)  \right)  \operatorname{div}u\\
&  =-\left(  \overline{\tau:\nabla u}-\tau:\nabla u\right)  \leq0,
\end{align*}
but by multiplication with $H^{-1}\left(  \overline{H\left(  \rho\right)
}-H\left(  \rho\right)  \right)  $ or $p^{-1}\left(  \overline{H\left(
\rho\right)  }-H\left(  \rho\right)  \right)  $ the left hand-side cannot be
written in conservative form.
\end{remark}

\section{Weak stability of sequences of global weak
solutions\label{WeakStability}}

\subsection{Classical functional analysis tools\label{Section_tools}}

This section is devoted to a quick recall of the main results from functional
analysis that we need in order to justify the computations done above. First,
we introduce a new function
\begin{equation}
g_{\varepsilon}=g\ast\omega_{\varepsilon}(x)\qquad\hbox{ with }\qquad
\omega_{\varepsilon}=\frac{1}{\varepsilon^{d}}\omega(\frac{x}{\varepsilon})
\label{notation_approx}%
\end{equation}
with $\omega$ a smooth nonnegative even function compactly supported in the
space ball of radius $1$ and with integral equal to 1. We recall the following
classical analysis result
\[
\lim_{\varepsilon\rightarrow0}\left\Vert g_{\varepsilon}-g\right\Vert
_{L^{q}(0,T;L^{p}({\mathbb{T}}^{3}))}=0.
\]
Next let us recall the following comutator estimate which was obtaiend for the
first time by DiPerna and Lions:

\begin{proposition}
\label{Prop_ren1}Consider $\beta\in(1,\infty)$ and $\left(  a,b\right)  $ such
that $a\in L^{\beta}\left(  \left(  0,T\right)  \times\mathbb{T}^{3}\right)  $
and $b,\nabla b\in L^{p}\left(  \left(  0,T\right)  \times\mathbb{T}%
^{3}\right)  $ where $\frac{1}{s}=\frac{1}{\beta}+\frac{1}{p}\leq1$. Then, we
have%
\[
\lim r_{\varepsilon}\left(  a,b\right)  =0\text{ in }L^{s}\left(  \left(
0,T\right)  \times\mathbb{T}^{3}\right)
\]
where
\begin{equation}
r_{\varepsilon}\left(  a,b\right)  =\partial_{i}\left(  a_{\varepsilon
}b\right)  -\partial_{i}\left(  \left(  ab\right)  _{\varepsilon}\right)  .
\label{def_reminder}%
\end{equation}

\end{proposition}

Whenever we have a \textit{regular solution} for the transport equation
\begin{equation}
\partial_{t}\rho+\operatorname{div}\left(  \rho u\right)  =0,
\label{transport_eq}%
\end{equation}
then, multiplying the former equation with $b^{\prime}\left(  \rho\right)  $
gives
\begin{equation}
\partial_{t}b\left(  \rho\right)  +\operatorname{div}\left(  b\left(
\rho\right)  u\right)  +\left\{  \rho b^{\prime}\left(  \rho\right)  -b\left(
\rho\right)  \right\}  \operatorname{div}u=0. \label{renorm}%
\end{equation}
The following proposition gives us a framework for justifying this
computations where $\rho$ is just a Lebesgue function.

\begin{proposition}
\label{Prop_ren2}Consider $2\leq\beta<\infty$ and $\lambda_{0},\lambda_{1}$
such that $\lambda_{0}<1$ and $-1\leq\lambda_{1}\leq\beta/2-1$. Also, consider
$\rho\in L^{\beta}\left(  \left(  0,T\right)  \times\mathbb{T}^{3}\right)  $,
$\rho\geq0$ a.e. and $u,\nabla u\in L^{2}\left(  \left(  0,T\right)
\times\mathbb{T}^{3}\right)  $ verifying the transport equation $\left(
\text{\ref{transport_eq}}\right)  $ in the sense of distributions. Then, for
any function $b\in C^{0}\left(  [0,\infty)\right)  \cap C^{1}\left(  \left(
0,\infty\right)  \right)  $ such that%
\[
\left\{
\begin{array}
[c]{l}%
b^{\prime}\left(  t\right)  \leq ct^{-\lambda_{0}}\text{ for }t\in(0,1],\\
\left\vert b^{\prime}\left(  t\right)  \right\vert \leq ct^{\lambda_{1}}\text{
for }t\geq1
\end{array}
\right.
\]
Then, equation $\left(  \text{\ref{renorm}}\right)  $ holds in the sense of distributions.
\end{proposition}

\noindent The proof of the above results follow by adapting in a
straightforward manner lemmas $6.7.$ and $6.9$ from the book of
Novotny-Stra\v{s}kraba \cite{NoSt} pages $304-308$.

\subsection{Estimates for bounded-energy weak
solutions\label{Aprioriestimates}}

Let us begin this section by recalling the basic a priori estimates for
(regular) solutions for the system%
\begin{equation}
\left\{
\begin{array}
[c]{l}%
\partial_{t}\rho+\operatorname{div}\left(  \rho u\right)  =0,\\
-\mathrm{div}\,\tau+\nabla\rho^{\gamma}=f,\\
\rho_{|t=0}=\rho_{0},
\end{array}
\right.  \label{NSC0}%
\end{equation}
with $\tau_{ij}=A_{ijkl}(t,x)D_{kl}(u)$ and
\[
\int_{\mathbb{T}^{3}}u\left(  t\right)  =\int_{\mathbb{T}^{3}}f\left(
t\right)  =0.
\]

Observe that we have set the adiabatic constant $a$ to equal to one just for
the sake of simplicity in the computations that follow.

First, of course, we have the mass conservation:%
\begin{equation}
\int_{\mathbb{T}^{3}}\rho(t)=\int_{\mathbb{T}^{3}}\rho|_{t=0}=\int
_{\mathbb{T}^{3}}\rho_{0}, \label{mass}%
\end{equation}
for all $t>0$ which follows by integrating the first equation of $\left(
\text{\ref{NSC0}}\right)  $. Next, by multiplying the velocity equation with
$u$ and integrating in space and time we get that%
\begin{align}
\int_{\mathbb{T}^{3}}\rho^{\gamma}\left(  t\right)  +\int_{0}^{t}%
\int_{\mathbb{T}^{3}}\tau &  :\nabla u\leq\int_{\mathbb{T}^{3}}\rho
_{0}^{\gamma}+\int_{0}^{t}\int_{\mathbb{T}^{3}}uf\label{energy}\\
&  \leq\int_{\mathbb{T}^{3}}\rho_{0}^{\gamma}+\left\Vert u\right\Vert
_{L_{t}^{2}L^{6}}\left\Vert f\right\Vert _{L_{t}^{2}L^{\frac{6}{5}}}.
\label{energy_222}%
\end{align}
The coercivity hypothesis $\left(  \text{\ref{H4}}\right)  $%
\[
c\int_{\mathbb{T}^{3}}|\nabla u|^{2}\leq\int_{\mathbb{T}^{3}}\tau:\nabla u,
\]
with $c>0$, the zero mean value on $u$, the K\"{o}rn inequality and Sobolev
embedding allows us to conclude that
\[
\rho\in L^{\infty}\left(  0,T;L^{\gamma}\left(  \mathbb{T}^{3}\right)
\right)  \text{, }\qquad u\in L^{2}(0,T;H^{1}\left(  \mathbb{T}^{3}\right)  )
\]
with%
\begin{equation}
\int_{\mathbb{T}^{3}}\rho^{\gamma}\left(  t\right)  +\int_{0}^{t}%
\int_{\mathbb{T}^{3}}|\nabla u|^{2}\leq C\left(  c\right)  \left(  \left\Vert
\rho_{0}\right\Vert _{L^{\gamma}}^{\gamma}+\int_{0}^{t}\left\Vert f\left(
\tau\right)  \right\Vert _{L^{\frac{6}{5}}}^{2}d\tau\right)  ,
\label{energy_2}%
\end{equation}
for all $t\geq0$ where $C\left(  c\right)  $ is a constant depending only on
the coercivity constant appearing in $\left(  \text{\ref{H4}}\right)  $.

Of course, the previous computations hold for regular solutions. It is to be
expected however that any reasonably physical solution to $\left(
\text{\ref{NSC0}}\right)  $ would verify the mass conservation and the energy
inequality. Thus, we introduce the following

\begin{definition}
Consider $f\in L^{2}\left(  0,T;L^{\frac{6}{5}}(\mathbb{T}^{3})\right)  $. A
pair $$\left(  \rho,u\right)  \in L^{\infty}\left(  0,T;L^{\gamma}\left(
\mathbb{T}^{3}\right)  \right)  \cap\mathcal{C}([0,T];L_{weak}^{\gamma
}(\mathbb{T}^{3}))\times L^{2}\left(  0,T;H^{1}\left(  \mathbb{T}^{3}\right)
\right)  $$ is called a bounded energy weak-solution for $\left(
\text{\ref{NSC0}}\right)  $ if it is a solution in the sense of distributions
for $\left(  \text{\ref{NSC0}}\right)  $ which moreover verifies the mass
conservation identity $\left(  \text{\ref{mass}}\right)  $ along with the
energy inequality $\left(  \text{\ref{energy}}\right)  $.
\end{definition}

This definition of bounded energy weak-solutions is consistent with the one we
find in Novotny-Stra\v{s}kraba \cite{NoSt} page $316.$

Of course, a bounded energy weak-solution for $\left(  \text{\ref{NSC0}%
}\right)  $ also verifies $\left(  \text{\ref{energy_2}}\right)  $. It turns
out that bounded energy weak-solutions verify some extra integrability
properties. More precisely, we have

\begin{proposition}
\label{ExtraEstimates} Consider $\left(  \rho,u\right)  \in L^{\infty}\left(
0,T;L^{\gamma}\left(  \mathbb{T}^{3}\right)  \right)  \times L^{2}\left(
0,T;H^{1}\left(  \mathbb{T}^{3}\right)  \right)  $ a bounded energy
weak-solution for $\left(  \text{\ref{NSC0}}\right)  $. Then, we have that%
\begin{equation}
\left.
\begin{array}
[c]{l}%
\left\Vert \rho^{\gamma}\right\Vert _{L_{t,x}^{2}}\leq C\left(  c,\gamma
\right)  \left(  \sqrt{t}+\max\left\{  1,\left\Vert A\right\Vert _{L^{\infty}%
}\right\}  \right)  \left(  \left\Vert \rho_{0}\right\Vert _{L^{\gamma}%
}^{\frac{\gamma}{2}}+\left\Vert f\right\Vert _{L_{t}^{2}L^{\frac{6}{5}}%
}\right)  ,\\
\left\Vert \partial_{t}u\right\Vert _{L^{1}(0,T;L^{\frac{3}{2}-\delta
}({\mathbb{T}}^{3}))}\leq C\left(  c,\gamma\right)  \left(  \sqrt{t}%
+\max\left\{  1,\left\Vert A\right\Vert _{L_{t,x}^{\infty}}\right\}  \right)
\left(  \left\Vert \rho_{0}\right\Vert _{L^{\gamma}}^{\gamma}+\left\Vert
f\right\Vert _{L_{t}^{2}L^{\frac{6}{5}}}^{2}\right) \\
\text{ \ \ \ \ \ \ \ \ \ \ \ \ \ \ \ \ \ \ \ \ \ \ \ \ \ \ \ \ \ \ \ \ \ }%
+C\left(  c,\gamma\right)  \sqrt{t}\left(  1+\left\Vert \partial
_{t}A\right\Vert _{L_{t,x}^{\infty}}\right)  \left(  \left\Vert \rho
_{0}\right\Vert _{L^{\gamma}}^{\frac{\gamma}{2}}+\left\Vert \left(
f,\partial_{t}f\right)  \right\Vert _{L_{t}^{2}L^{\frac{6}{5}}}\right)  ,
\end{array}
\right.  \label{rho_gamma_extra_and_derivative_u}%
\end{equation}
where $C\left(  c,\gamma\right)  $ depends only on $c$ and $\gamma$ and
$\delta\in\left(  0,1/2\right)  $ is the constant appearing in
\textrm{$\left(  \text{\ref{H4}}\right)  $}.
\end{proposition}

\subparagraph{Proof or Proposition \ref{ExtraEstimates}:}

The integrability assumptions for the weak solution $\left(  \rho,u\right)  $
ensure that for all $\psi\in\left[  L^{2}\left(  0,T;H^{1}\left(
\mathbb{T}^{3}\right)  \right)  \right]  ^{3}$ we have that%
\[
\int_{0}^{t}%
{\displaystyle\int_{\mathbb{T}^{3}}}
\rho^{\gamma}\operatorname{div}\psi=\int_{0}^{t}%
{\displaystyle\int_{\mathbb{T}^{3}}}
\tau:\nabla\psi+\int_{0}^{t}%
{\displaystyle\int_{\mathbb{T}^{3}}}
f\psi
\]

Taking $\phi\in L^{2}\left(  (0,T)\times\mathbb{T}^{3}\right)  $ and
considering a test function $\psi$ such that%
\[
\Delta\psi=\nabla\phi\text{ with}%
{\displaystyle\int_{\mathbb{T}^{3}}}
\psi=0,
\]
we get that
\[
\operatorname{div}\psi=\phi-%
{\displaystyle\int_{\mathbb{T}^{3}}}
\phi,
\]
and owing to $A(t,x)\in W^{1,\infty}((0,T)\times{\mathbb{T}}^{3}))^{3\times3}$
along with the energy estimate $\left(  \text{\ref{energy_2}}\right)  $, we
get that
\begin{align*}
\int_{0}^{t}%
{\displaystyle\int_{\mathbb{T}^{3}}}
\rho^{\gamma}\phi &  =\int_{0}^{t}%
{\displaystyle\int_{\mathbb{T}^{3}}}
\phi%
{\displaystyle\int_{\mathbb{T}^{3}}}
\rho^{\gamma}+\int_{0}^{t}%
{\displaystyle\int_{\mathbb{T}^{3}}}
\tau:\nabla\psi+\int_{0}^{t}%
{\displaystyle\int_{\mathbb{T}^{3}}}
f\psi\\
&  \leq C\left(  c,\gamma\right)  \left(  \sqrt{t}+\max\left\{  1,\left\Vert
A\right\Vert _{L^{\infty}}\right\}  \right)  \left(  \left\Vert \rho
_{0}\right\Vert _{L^{\gamma}}^{\frac{\gamma}{2}}+\left\Vert f\right\Vert
_{L_{t}^{2}L^{\frac{6}{5}}}\right)  \left\Vert \phi\right\Vert _{L_{t,x}^{2}}%
\end{align*}
and thus we get that%
\begin{equation}
\rho^{\gamma}\in L^{2}\left(  \left(  0,T\right)  \times\mathbb{T}^{3}\right)
, \label{extra_integrability}%
\end{equation}
verifying uniform bound announced in the first relation of $\left(
\text{\ref{rho_gamma_extra_and_derivative_u}}\right)  $.

We prove now the estimate for the time derivative of $\partial_{t}u$. We can
recover time regularity for $u$ by proceeding in the following way. We write
that%
\begin{align*}
-{\mathcal{A}}\partial_{t}u  &  =\mathrm{div}(\partial_{t}A(t,x)D(u))+\partial
_{t}f-\nabla\partial_{t}\rho^{\gamma}\\
&  =\mathrm{div}(\partial_{t}A(t,x)D(u))+\partial_{t}f\\
&  \hskip.5cm+\nabla\operatorname{div}\left(  \rho^{\gamma}u-%
{\displaystyle\int_{\mathbb{T}^{3}}}
\rho^{\gamma}u\right)  +\left(  \gamma-1\right)  \nabla\left(  \rho^{\gamma
}\operatorname{div}u-%
{\displaystyle\int_{\mathbb{T}^{3}}}
\rho^{\gamma}\operatorname{div}u\right)  .
\end{align*}
where the passage from the second line to the third is justified by
Proposition $\left(  \text{\ref{Prop_ren2}}\right)  $ which of course, can be
applied owing to the fact that we recover $\left(
\text{\ref{extra_integrability}}\right)  $. Above, the first two terms behave
better and thus taking advandtage of the linearity of the operator
$-{\mathcal{A}}$ it is more convenient to separate $\partial_{t}u$ in two
parts and estimate them separetly. To this end, consider $\phi$ with
$\int_{\mathbb{T}^{3}}\phi=0$, such that%
\[
-\mathcal{A}\phi{\mathcal{=}}\mathrm{div}(\partial_{t}A(t,x)D(u))+\partial
_{t}f
\]
Multiplying by $\phi$ we get that%
\begin{align*}
c\int_{0}^{t}\int_{\mathbb{T}^{3}}\left\vert \nabla\phi\right\vert ^{2}  &
\leq-\int_{0}^{t}\int_{\mathbb{T}^{3}}\phi\mathcal{A}\phi=-\int_{0}^{t}%
\int_{\mathbb{T}^{3}}\partial_{t}A(t,x)D(u)\nabla\phi+\int_{0}^{t}%
\int_{\mathbb{T}^{3}}\partial_{t}f\phi\\
&  \leq\frac{1}{8c}\int_{0}^{t}\int_{\mathbb{T}^{3}}\left\vert \partial
_{t}A(t,x)D(u)\right\vert ^{2}+\frac{C^{2}}{8c}\int_{0}^{t}\left\Vert
\partial_{t}f\right\Vert _{L^{\frac{6}{5}}}^{2}+\frac{c}{2}\int_{0}^{t}%
\int_{\mathbb{T}^{3}}\left\vert \nabla\phi\right\vert ^{2}%
\end{align*}
where $C$ is the constant appearing in the Sobolev inequality and thus, we get
that%
\begin{equation}
\frac{c}{2}\int_{0}^{t}\int_{\mathbb{T}^{3}}\left\vert \nabla\phi\right\vert
^{2}\leq\frac{1}{8c}\int_{0}^{t}\int_{\mathbb{T}^{3}}\left\vert \partial
_{t}A(t,x)D(u)\right\vert ^{2}+\frac{C^{2}}{8c}\int_{0}^{t}\left\Vert
\partial_{t}f\right\Vert _{L^{\frac{6}{5}}}^{2}. \label{bound on phi}%
\end{equation}
It remains to estimate $\partial_{t}u-\phi$ which verifies%
\[
{\mathcal{A}}\left(  \partial_{t}u-\phi\right)  =-\nabla\operatorname{div}%
\left(  \rho^{\gamma}u-\int\rho^{\gamma}u\right)  -\left(  \gamma-1\right)
\nabla\left(  \rho^{\gamma}\operatorname{div}u-\int\rho^{\gamma}%
\operatorname{div}u\right)  .
\]
We will use a periodic variant of the following result due to Stampacchia and
for more general second order elliptic equation to Boccardo-Gallou\"{e}t that
can be found for instance in \cite{Ponce} Proposition $5.1.$ page $77$. Let
$\psi$ be the solution of
\[
-\Delta\psi=f\text{ with }\psi_{|\partial\Omega}=0,
\]
where $f\in L^{1}(\Omega)$ with $\Omega$ a smooth bounded domain then we have
that%
\begin{equation}
\left\Vert \nabla\psi\right\Vert _{L^{r}(\Omega)}\leq C_{\delta}\left\Vert
f\right\Vert _{L^{1}(\Omega)} \label{elliptique_Ponce}%
\end{equation}
for all $r\in\lbrack1,3/2)$. The periodic version reads as follows: let $\psi$
a solution of
\[
-\Delta\psi=f\hbox{ with }f\in L^{1}({\mathbb{T}}^{3})\hbox{ and }\int
_{{\mathbb{T}}^{3}}f=0
\]
then \eqref{elliptique_Ponce} is satisfied, see Theorem
\ref{Stampacchia_Gallouet} from the Appendix for a proof. As $\rho^{\gamma
}\operatorname{div}u\in L^{1}((0,T)\times{\mathbb{T}}^{3})$, let us consider
$\psi$ the solution of
\[
-\Delta\psi\left(  \rho,u\right)  =\rho^{\gamma}\operatorname{div}%
u-\int_{\mathbb{T}^{3}}\rho^{\gamma}\operatorname{div}u
\]
which verifies that%
\[
\left\Vert \nabla\psi\left(  \rho,u\right)  \right\Vert _{L^{1}(0,T:L^{\frac
{3}{2}-\delta}({\mathbb{T}}^{3}))}\leq C_{\delta}\left\Vert \rho^{\gamma
}\operatorname{div}u\right\Vert _{L^{1}(0,T;L^{1}({\mathbb{T}}^{3}))}\leq
C_{\delta}\left\Vert \rho^{\gamma}\right\Vert _{L^{2}((0,T)\times{\mathbb{T}%
}^{3})}\left\Vert \operatorname{div}u\right\Vert _{L^{2}((0,T)\times
{\mathbb{T}}^{3})}.
\]
where $\delta\in\left(  0,1/2\right)  $ is the constant appearing in $\left(
\text{\ref{H4}}\right)  $. But then, we may write that%
\begin{align*}
{\mathcal{A}}\left(  \partial_{t}u-\phi\right)   &  =-\nabla\operatorname{div}%
\left(  \rho^{\gamma}u\right)  -\left(  \gamma-1\right)  \nabla\left(
\rho^{\gamma}\operatorname{div}u\right) \\
&  =\nabla\operatorname{div}\left(  \rho^{\gamma}u\right)  +\left(
\gamma-1\right)  \nabla\operatorname{div}\nabla\psi\left(  \rho,u\right)
\end{align*}
and using hypothesis $\left(  \text{\ref{H4}}\right)  $ we arrive at
\begin{align}
\left\Vert \left(  \partial_{t}u-\phi\right)  \right\Vert _{L^{1}%
(0,T;L^{\frac{3}{2}-\delta}({\mathbb{T}}^{3}))}  &  \leq\left\Vert
\rho^{\gamma}u-\int_{\mathbb{T}^{3}}\rho^{\gamma}u\right\Vert _{L^{1}%
(0,T;L^{\frac{3}{2}}({\mathbb{T}}^{3}))}+\left\Vert \nabla\psi\left(
\rho,u\right)  \right\Vert _{L^{1}(0,T;L^{\frac{3}{2}-\delta}({\mathbb{T}}%
^{3}))}\nonumber\\
&  \leq\left\Vert \rho^{\gamma}\right\Vert _{L^{2}((0,T)\times{\mathbb{T}}%
^{3})}\left\Vert u\right\Vert _{L^{2}(0,T;L^{6}({\mathbb{T}}^{3}))}+\left\Vert
\rho^{\gamma}\right\Vert _{L^{2}((0,T)\times{\mathbb{T}}^{3})}\left\Vert
\operatorname{div}u\right\Vert _{L^{2}((0,T)\times{\mathbb{T}}^{3}%
)}\nonumber\\
&  \leq\left\Vert \rho^{\gamma}\right\Vert _{L^{2}((0,T)\times{\mathbb{T}}%
^{3})}\left\Vert \nabla u\right\Vert _{L^{2}((0,T)\times{\mathbb{T}}^{3})}.
\label{bound on remainder}%
\end{align}
We get a uniform bound for $\partial_{t}u$ in $L^{1}\left(  0,T;L^{3/2-}%
(\mathbb{T}^{3}\right)  )$ by combining estimates $\left(
\text{\ref{bound on phi}}\right)  $ and $\left(
\text{\ref{bound on remainder}}\right)  $ in the following manner%
\begin{align*}
\left\Vert \partial_{t}u\right\Vert _{L^{1}(0,T;L^{\frac{3}{2}-\delta
}({\mathbb{T}}^{3}))}  &  \leq\left\Vert \left(  \partial_{t}u-\phi\right)
\right\Vert _{L^{1}(0,T;L^{\frac{3}{2}-\delta}({\mathbb{T}}^{3}))}+\left\Vert
\phi\right\Vert _{L^{1}(0,T;L^{\frac{3}{2}-\delta}({\mathbb{T}}^{3}))}\\
&  \leq\left\Vert \rho^{\gamma}\right\Vert _{L^{2}((0,T)\times{\mathbb{T}}%
^{3})}\left\Vert \nabla u\right\Vert _{L^{2}((0,T)\times{\mathbb{T}}^{3}%
)}+\sqrt{t}\left\Vert \phi\right\Vert _{L^{2}(0,T;L^{6}({\mathbb{T}}^{3}))}\\
&  \leq\left\Vert \rho^{\gamma}\right\Vert _{L^{2}((0,T)\times{\mathbb{T}}%
^{3})}\left\Vert \nabla u\right\Vert _{L^{2}((0,T)\times{\mathbb{T}}^{3}%
)}+\sqrt{t}\left\Vert \nabla\phi\right\Vert _{L^{2}((0,T)\times{\mathbb{T}%
}^{3})}\\
&  \leq C\left(  c,\gamma\right)  \left(  \sqrt{t}+\max\left\{  1,\left\Vert
A\right\Vert _{L_{t,x}^{\infty}}\right\}  \right)  \left(  \left\Vert \rho
_{0}\right\Vert _{L^{\gamma}}^{\gamma}+\left\Vert f\right\Vert _{L_{t}%
^{2}L^{\frac{6}{5}}}^{2}\right) \\
&  +C\left(  c,\gamma\right)  \sqrt{t}\left(  1+\left\Vert \partial
_{t}A\right\Vert _{L_{t,x}^{\infty}}\right)  \left(  \left\Vert \rho
_{0}\right\Vert _{L^{\gamma}}^{\frac{\gamma}{2}}+\left\Vert \left(
f,\partial_{t}f\right)  \right\Vert _{L_{t}^{2}L^{\frac{6}{5}}}\right)
\end{align*}
which is exactly the estimate $\left(
\text{\ref{rho_gamma_extra_and_derivative_u}}\right)  $. Of course combining
this information with the energy inequality $\left(  \text{\ref{energy}%
}\right)  $ we obtain an uniform bound for
\[
u\in L^{2}(0,T;H^{1}({\mathbb{T}}^{3}))\cap W^{1,1}(0,T;L^{3/2-\delta
}({\mathbb{T}}^{3})).
\]
This ends the proof of Proposition \ref{ExtraEstimates}.

\begin{remark}
Also, for later purposes it is convenient to observe that we actually proved
that if%
\begin{equation}
-{\mathcal{A}}u=\operatorname{div}F \label{Au=div}%
\end{equation}
then Hypothesis $\left(  \text{\ref{H4}}\right)  $ made on the operator
${\mathcal{A}}$ implies that there exists some constant $C$ such that%
\begin{equation}
\left\Vert \nabla u\right\Vert _{L^{\frac{3}{2}-\delta}({\mathbb{T}}^{3})}\leq
C\left\Vert F\right\Vert _{L^{1}({\mathbb{T}}^{3})}. \label{EstimateL1}%
\end{equation}
for any $u,F$ verifying $\left(  \text{\ref{Au=div}}\right)  $.
\end{remark}

\begin{remark}
The previous estimates are not all available in the case of the full
compressible Navier-Stokes system. For instance we do not have control on the
time derivative of the velocity and $\rho^{\gamma}$ is not square integrable:
we control only $\partial_{t}(\rho u)$ in $L^{1}(0,T;H^{-1}({\mathbb{T}}%
^{3}))$ allowing to get compactness on $\sqrt{\rho}u$ in $L^{2}((0,T)\times
{\mathbb{T}}^{3}))$ and we gain extra integrability $\rho^{\gamma+\theta}\in
L^{1}((0,T)\times{\mathbb{T}}^{3})$ for $0<\theta<2\gamma/3-1$.
\end{remark}

\subsection{Weak stability of solutions of $\left(  \text{\ref{equations}%
}\right)  $\label{Weakstability}}

\label{weakstab} The aim of this section is to provide the arguments that
render rigorous the formal computations presented in Section \ref{Newapproach}%
. Let us temporarily include an extra potential source term in the system:%

\begin{equation}
\left\{
\begin{array}
[c]{l}%
\partial_{t}\rho+\operatorname{div}\left(  \rho u\right)  =0,\\
-\mathrm{div}\,\tau+\nabla\rho^{\gamma}=\nabla g+f.
\end{array}
\right.  \label{NSC1}%
\end{equation}
As we saw in Section \ref{Section_tools} under certain integrability
conditions one may conclude that $\rho^{\gamma}$ verifies the following
equation :%
\[
\partial_{t}\rho^{\gamma}+\operatorname{div}\left(  \rho^{\gamma}u\right)
+\left(  \gamma-1\right)  \rho^{\gamma}\operatorname{div}u=0.
\]
Of course, the result of Proposition \ref{Prop_ren2} that allows us to write
the above equation does not take in account the structure of the system
$\left(  \text{\ref{NSC1}}\right)  $. In the following, we propose a more
accurate result taking in consideration the equation of the velocity.

\begin{proposition}
\label{Prop_ren_new}\bigskip Consider $f\in L^{2}(0,T;L^{\frac{6}{5}%
}(\mathbb{T}^{3}))$, $g\in L^{2}\left(  \left(  0,T\right)  \times
\mathbb{T}^{3}\right)  $ and $\left(  \rho,u\right)  $ a bounded energy
weak-solution of \eqref{NSC1}. Then, one has that%
\begin{equation}
\tfrac{1}{\gamma-1}\left\{  \partial_{t}\rho^{\gamma}+\gamma\operatorname{div}%
\left(  \rho^{\gamma}u\right)  \right\}  =\mathrm{div}(\tau:u)-\tau:\nabla
u+uf+\operatorname{div}\left(  ug\right)  -g\operatorname{div}u.
\label{Prop_ren_new_1}%
\end{equation}
in the sense of distributions.
\end{proposition}

\begin{remark}
In order to prove Proposition \ref{Prop_ren_new} we do not require regularity
on the time derivative of $f$ as it is needed in order to obtain the a priori
estimates for $\partial_{t}u$, see Proposition \ref{ExtraEstimates}.
\end{remark}

\begin{remark}
Proposition \ref{Prop_ren_new} is valid for all tensor fields $\tau\in
L^{2}\left(  \left(  0,T\right)  \times\mathbb{T}^{3}\right)  $.
\end{remark}

\subparagraph{Proof of \ref{Prop_ren_new}:}

The proof uses the regularizing the techniques introduced by Lions in
\cite{Li}, see also the book of Novotny and Stra\v{s}kraba (\cite{NoSt}).
Recall the notation introduced in $\left(  \text{\ref{notation_approx}%
}\right)  $ and $\left(  \text{\ref{def_reminder}}\right)  $ and let us write
\[
\partial_{t}\rho_{\varepsilon}+\operatorname{div}\left(  \rho_{\varepsilon
}u\right)  =r_{\varepsilon}\left(  \rho,u\right)
\]
which by multiplying with $\gamma(\rho_{\varepsilon})^{\gamma-1}$ yields%
\[
\partial_{t}\left(  \rho_{\varepsilon}\right)  ^{\gamma}+\operatorname{div}%
\left(  \left(  \rho_{\varepsilon}\right)  ^{\gamma}u\right)  +(\gamma
-1)\left(  \rho_{\varepsilon}\right)  ^{\gamma}\operatorname{div}%
u=\gamma\,r_{\varepsilon}\left(  \rho,u\right)  \left(  \rho_{\varepsilon
}\right)  ^{\gamma-1}.
\]
Let us rewrite the above equation in the following manner:
\begin{gather*}
\partial_{t}\left(  \rho_{\varepsilon}\right)  ^{\gamma}+\operatorname{div}%
\left(  \left(  \rho_{\varepsilon}\right)  ^{\gamma}u\right)  +(\gamma
-1)\left\{  \left(  \rho_{\varepsilon}\right)  ^{\gamma}-(\rho^{\gamma
})_{\varepsilon^{\prime}}\right\}  \operatorname{div}u+(\gamma-1)(\rho
^{\gamma})_{\varepsilon^{\prime}}\left\{  \operatorname{div}%
u-\operatorname{div}u_{\varepsilon^{\prime}}\right\} \\
+(\gamma-1)(\rho^{\gamma})_{\varepsilon^{\prime}}\operatorname{div}%
u_{\varepsilon^{\prime}}=\gamma r_{\varepsilon}\,\left(  \rho,u\right)
\left(  \rho_{\varepsilon}\right)  ^{\gamma-1}.
\end{gather*}
Next, we observe that owing to the second equation of $\left(
\text{\ref{NSC1}}\right)  $ we get that%
\begin{align*}
(\rho^{\gamma})_{\varepsilon^{\prime}}\operatorname{div}u_{\varepsilon
^{\prime}}  &  =\operatorname{div}\left(  (\rho^{\gamma})_{\varepsilon
^{\prime}}u_{\varepsilon^{\prime}}\right)  -u_{\varepsilon^{\prime}}%
\nabla(\rho^{\gamma})_{\varepsilon^{\prime}}\\
&  =\operatorname{div}\left(  (\rho^{\gamma})_{\varepsilon^{\prime}%
}u_{\varepsilon^{\prime}}\right)  -u_{\varepsilon^{\prime}}\mathrm{div}%
\tau_{\varepsilon^{\prime}}-u_{\varepsilon^{\prime}}\nabla g_{\varepsilon
^{\prime}}-u_{\varepsilon^{\prime}}f_{\varepsilon^{\prime}}\\
&  =\operatorname{div}\left(  (\rho^{\gamma})_{\varepsilon^{\prime}%
}u_{\varepsilon^{\prime}}\right)  -\mathrm{div}(\tau_{\varepsilon^{\prime}%
}:u_{\varepsilon^{\prime}})+\tau_{\varepsilon^{\prime}}:\nabla u_{\varepsilon
^{\prime}}-\operatorname{div}\left(  u_{\varepsilon^{\prime}}g_{\varepsilon
^{\prime}}\right)  +g_{\varepsilon^{\prime}}\operatorname{div}u_{\varepsilon
^{\prime}}-u_{\varepsilon^{\prime}}f_{\varepsilon^{\prime}}%
\end{align*}
and thus, we may write that%
\begin{gather*}
\tfrac{1}{\gamma-1}\left\{  \partial_{t}\left(  \rho_{\varepsilon}\right)
^{\gamma}+\operatorname{div}\left(  \left(  \rho_{\varepsilon}\right)
^{\gamma}u\right)  \right\}  +\left\{  \left(  \rho_{\varepsilon}\right)
^{\gamma}-(\rho^{\gamma})_{\varepsilon^{\prime}}\right\}  \operatorname{div}%
u+(\rho^{\gamma})_{\varepsilon^{\prime}}\left\{  \operatorname{div}%
u-\operatorname{div}u_{\varepsilon^{\prime}}\right\} \\
+\operatorname{div}\left(  (\rho^{\gamma})_{\varepsilon^{\prime}%
}u_{\varepsilon^{\prime}}\right)  -\mathrm{div}(\tau_{\varepsilon^{\prime}%
}:u_{\varepsilon^{\prime}})+\tau_{\varepsilon^{\prime}}:\nabla u_{\varepsilon
^{\prime}}-\operatorname{div}\left(  u_{\varepsilon^{\prime}}g_{\varepsilon
^{\prime}}\right)  +g_{\varepsilon^{\prime}}\operatorname{div}u_{\varepsilon
^{\prime}}-u_{\varepsilon^{\prime}}f_{\varepsilon^{\prime}}\\
=\tfrac{\gamma}{\gamma-1}r_{\varepsilon}\left(  \rho,u\right)  \left(
\rho_{\varepsilon}\right)  ^{\gamma-1}.
\end{gather*}
Using the strong convergence properties of the convolution, Proposition
\ref{Prop_ren1} along with the fact that bounded energy weak-solutions also
satisfy $\rho\in L^{2\gamma}\left(  \left(  0,T\right)  \times\mathbb{T}%
^{3}\right)  $ we get that
\[
\left\{
\begin{array}
[c]{l}%
\left(  \rho_{\varepsilon}\right)  ^{\gamma}\rightarrow\rho^{\gamma}\text{ in
}L^{2}\left(  (0,T)\times\mathbb{T}^{3}\right)  \text{ for }\varepsilon
\rightarrow0,\\
\left(  \rho_{\varepsilon}\right)  ^{\gamma}u\rightarrow\rho^{\gamma}u\text{
in }L^{1}\left(  (0,T)\times\mathbb{T}^{3}\right)  \text{ for }\varepsilon
\rightarrow0,\\
(\rho^{\gamma})_{\varepsilon^{\prime}}\left\{  \operatorname{div}%
u-\operatorname{div}u_{\varepsilon^{\prime}}\right\}  \rightarrow0\text{ in
}L^{1}\left(  (0,T)\times\mathbb{T}^{3}\right)  \text{ for }\varepsilon
^{\prime}\rightarrow0\\
(\rho^{\gamma})_{\varepsilon^{\prime}}\operatorname{div}u_{\varepsilon
^{\prime}}\rightarrow\rho^{\gamma}\operatorname{div}u\text{ in }L^{1}\left(
(0,T)\times\mathbb{T}^{3}\right)  \text{ for }\varepsilon^{\prime}%
\rightarrow0,\\
\tau_{\varepsilon^{\prime}}:u_{\varepsilon^{\prime}}\rightarrow\tau:u\text{
and }\tau_{\varepsilon^{\prime}}:\nabla u_{\varepsilon^{\prime}}%
\rightarrow\tau:u\text{ in }L^{1}\left(  (0,T)\times\mathbb{T}^{3}\right)
\text{ for }\varepsilon^{\prime}\rightarrow0,\\
u_{\varepsilon^{\prime}}f_{\varepsilon^{\prime}}\rightarrow uf\text{ in }%
L^{1}\left(  (0,T)\times\mathbb{T}^{3}\right)  \text{ for }\varepsilon
^{\prime}\rightarrow0,\\
u_{\varepsilon^{\prime}}g_{\varepsilon^{\prime}}\rightarrow ug\text{ in }%
L^{1}\left(  (0,T)\times\mathbb{T}^{3}\right)  \text{ for }\varepsilon
^{\prime}\rightarrow0,\\
g_{\varepsilon^{\prime}}\operatorname{div}u_{\varepsilon^{\prime}}\rightarrow
g\operatorname{div}u\text{ in }L^{1}\left(  (0,T)\times\mathbb{T}^{3}\right)
\text{ for }\varepsilon^{\prime}\rightarrow0,\\
r_{\varepsilon}\left(  \rho,u\right)  \left(  \rho_{\varepsilon}\right)
^{\gamma-1}\rightarrow0\text{ in }L^{1}\left(  (0,T)\times\mathbb{T}%
^{3}\right)  \text{ for }\varepsilon\rightarrow0.
\end{array}
\text{ }\right.
\]
Consequently, we get that
\[
\tfrac{1}{\gamma-1}\left\{  \partial_{t}\rho^{\gamma}+\gamma\operatorname{div}%
\left(  \rho^{\gamma}u\right)  \right\}  =\mathrm{div}(\tau u)-\tau:\nabla
u+fu+\operatorname{div}\left(  gu\right)  -g\operatorname{div}u.
\]
This ends the proof of Proposition \ref{Prop_ren_new}. Next, we investigate
the weak stability of a sequence of solutions of system $\left(
\text{\ref{NSC1}}\right)  $. Our main results reads

\begin{theorem}
\label{main}\bigskip Consider a sequence of bounded energy weak-solutions
$\left(  \rho^{\varepsilon},u^{\varepsilon}\right)  _{\varepsilon>0}$ for
$\left(  \text{\ref{NSC1}}\right)  $ with initial data $\left(  \rho
_{0}^{\varepsilon}\right)  _{\varepsilon>0}\subset L^{\gamma}\left(
\mathbb{T}^{3}\right)  $, i.e.
\begin{equation}
\left\{
\begin{array}
[c]{l}%
\partial_{t}\rho^{\varepsilon}+\operatorname{div}\left(  \rho^{\varepsilon
}u^{\varepsilon}\right)  =0,\\
-\mathrm{div}\,\tau^{\varepsilon}+\nabla(\rho^{\varepsilon})^{\gamma
}=f^{\varepsilon},\\
\rho_{|t=0}^{\varepsilon}=\rho_{0}^{\varepsilon},
\end{array}
\right.  \label{eps_system}%
\end{equation}
with
\[
\tau_{ij}^{\varepsilon}=A_{ijkl}^{\varepsilon}(t,x)D_{kl}(u^{\varepsilon}),
\]
where%
\begin{equation}
\left\{
\begin{array}
[c]{l}%
\rho_{0}^{\varepsilon}\rightarrow\rho_{0}\text{ in }L^{\gamma}\left(
\mathbb{T}^{3}\right)  ,\\
A^{\varepsilon}(t,x)\rightarrow A(t,x)\hbox{ in }W^{1,\infty}((0,T)\times
{\mathbb{T}}^{3}),\\
f^{\varepsilon}\rightarrow f\hbox{ in }L^{2}(0,T;L^{\frac{6}{5}}({\mathbb{T}%
}^{3}))).
\end{array}
\right.  \label{uniform}%
\end{equation}
Then, there exists $\left(  \rho,u\right)  \in$ $L^{2\gamma}\left(
(0,T)\times\mathbb{T}^{3}\right)  \times$ $\left[  L^{2}(0,T;H^{1}%
(\mathbb{T}^{3}))\right]  ^{3}$ such that modulo a subsequence we have%
\begin{equation}
\left\{
\begin{array}
[c]{l}%
\rho^{\varepsilon}\rightharpoonup\rho\text{ weakly in }L^{2\gamma}\left(
(0,T)\times\mathbb{T}^{3}\right)  ,\\
\rho^{\varepsilon}\rightarrow\rho\text{ in }L^{2\gamma-}\left(  (0,T)\times
\mathbb{T}^{3}\right)  ,\\
u^{\varepsilon}\rightharpoonup u\text{ weakly in }L^{2}(0,T;H^{1}%
(\mathbb{T}^{3}))\text{ }\\
u^{\varepsilon}\rightarrow u\text{ in }L^{2}((0,T)\times{\mathbb{T}}^{3})),
\end{array}
\right.  \label{convergence_of_sequences}%
\end{equation}
where $\left(  \rho,u\right)  $ verifies%
\begin{equation}
\left\{
\begin{array}
[c]{l}%
\partial_{t}\rho+\operatorname{div}\left(  \rho u\right)  =0,\\
-\mathrm{div}\tau+\nabla\rho^{\gamma}=\nabla f,\\
\rho_{|t=0}=\rho_{0}.
\end{array}
\right.  \label{limit_system}%
\end{equation}
with
\[
\tau_{ij}=A_{ijkl}(t,x)D_{kl}(u).
\]
Morever, the following energy bound holds a.e. $t\in\left(  0,T\right)  :$%
\begin{equation}
\int_{\mathbb{T}^{3}}\rho^{\gamma}\left(  t\right)  +\int_{0}^{t}%
\int_{\mathbb{T}^{3}}\tau:\nabla u\leq\int_{\mathbb{T}^{3}}\rho_{0}^{\gamma
}+\int_{0}^{t}\int_{\mathbb{T}^{3}}uf.
\end{equation}

\end{theorem}

\subparagraph{Proof of Theorem \ref{main}}

The information on the initial data $\left(
\text{\ref{convergence_of_sequences}}\right)  $ along with Proposition
\ref{extra_integrability} ensures that%

\[
\left\Vert \rho^{\varepsilon}\right\Vert _{L^{\infty}(0,T;L^{\gamma
}({\mathbb{T}}^{3}))\cap L^{2\gamma}((0,T)\times{\mathbb{T}}^{3})}+\left\Vert
u^{\varepsilon}\right\Vert _{L^{2}(0,T;H^{1}({\mathbb{T}}^{3}))\cap
W^{1,1}(0,T;L^{3/2-\delta}({\mathbb{T}}^{3}))}\leq C\left(  1+\sqrt{T}\right)
,
\]
for all $T>0$. The assumptions allow us to conclude that there exist three
functions $\left(  \rho,u\right)  $ and $\overline{\rho^{\gamma}}$ such that
up to a subsequence we have the following informations :
\begin{equation}
\left\{
\begin{array}
[c]{l}%
\rho^{\varepsilon}\rightharpoonup\rho\text{ weakly in }L^{2\gamma}\left(
(0,T)\times\mathbb{T}^{3}\right)  ,\\
\rho^{\varepsilon}\rightarrow\rho\text{ strongly in }\mathcal{C}%
([0,T];L_{weak}^{\gamma}(\mathbb{T}^{3}))\\
(\rho^{\varepsilon})^{\gamma}\rightharpoonup\overline{\rho^{\gamma}}\text{
weakly in }L^{2}\left(  (0,T)\times\mathbb{T}^{3}\right)  ,\\
\nabla u^{\varepsilon}\rightharpoonup\nabla u\text{ weakly in }L^{2}\left(
(0,T)\times\mathbb{T}^{3}\right)  ,\\
u^{\varepsilon}\rightarrow u\text{ strongly in }L^{2}\left(  (0,T)\times
\mathbb{T}^{3}\right)  .
\end{array}
\right.  \label{weak_conv1}%
\end{equation}
Moreover, we may take the above subsequence such as%
\begin{equation}
\left\{
\begin{array}
[c]{c}%
\tau^{\varepsilon}:\nabla u^{\varepsilon}\rightharpoonup\overline{\tau:\nabla
u}\text{ in }\mathcal{M}\left(  (0,T)\times\mathbb{T}^{3}\right)  \text{
and}\\
\tau:\nabla u\leq\overline{\tau:\nabla u}\text{ in the sense of measures}%
\end{array}
\right.  \label{weak_conv2}%
\end{equation}
using the weak lower semi-continuity of the viscous work: see hypothesis
$\left(  \text{\ref{H2}}\right)  $. All the above information allows us to
conclude that%
\begin{equation}
\left\{
\begin{array}
[c]{l}%
\partial_{t}\rho+\operatorname{div}\left(  \rho u\right)  =0,\\
-\mathrm{div}\,\tau+\nabla\overline{\rho^{\gamma}}=f,
\end{array}
\right.  \label{limit_not_yet_ident}%
\end{equation}
with
\[
\tau_{ij}=A_{ijkl}(t,x)D_{kl}(u).
\]
Of course, the most delicate part is to identify $\overline{\rho^{\gamma}}$
with $\rho^{\gamma}$. Let us observe that for any $\varepsilon>0$ , $\left(
\rho^{\varepsilon},u^{\varepsilon}\right)  $ verifies the hypothesis of
Proposition \ref{Prop_ren_new} and thus we infer that%
\begin{equation}
\tfrac{1}{\gamma-1}\left\{  \partial_{t}(\rho^{\varepsilon})^{\gamma}%
+\gamma\operatorname{div}\left(  (\rho^{\varepsilon})^{\gamma}u^{\varepsilon
}\right)  \right\}  =\mathrm{div}(\tau^{\varepsilon}:u^{\varepsilon}%
)-\tau^{\varepsilon}:\nabla u^{\varepsilon}+f^{\varepsilon}u^{\varepsilon}
\label{renorm_eps}%
\end{equation}
Moreover, using the information of relation $\left(  \text{\ref{weak_conv1}%
}\right)  $ we may pass to the limit in $\left(  \text{\ref{renorm_eps}%
}\right)  $ such as to obtain%
\begin{equation}
\tfrac{1}{\gamma-1}\left\{  \partial_{t}\overline{\rho^{\gamma}}%
+\gamma\operatorname{div}\left(  (\overline{\rho^{\gamma}}u\right)  \right\}
=\mathrm{div}(\tau:u)-\overline{\tau:\nabla u}+fu. \label{passage_limite}%
\end{equation}
Observing that we may put the system $\left(  \text{\ref{limit_not_yet_ident}%
}\right)  $ under the form
\begin{equation}
\left\{
\begin{array}
[c]{l}%
\partial_{t}\rho+\operatorname{div}\left(  \rho u\right)  =0,\\
-\mathrm{div}\,\tau+\nabla\rho^{\gamma}=\nabla(\rho^{\gamma}-\overline
{\rho^{\gamma}})+f
\end{array}
\right.
\end{equation}
with $\tau_{ij}=A_{ijkl}(t,x)D_{kl}(u)$ and using Proposition
\ref{Prop_ren_new} we write that%
\begin{align}
&  \tfrac{1}{\gamma-1}\left\{  \partial_{t}\rho^{\gamma}+\gamma
\operatorname{div}\left(  \rho^{\gamma}u\right)  \right\}  -\operatorname{div}%
\left(  u\left(  \rho^{\gamma}-\overline{\rho^{\gamma}}\right)  \right)
+\left(  \rho^{\gamma}-\overline{\rho^{\gamma}}\right)  \operatorname{div}%
u\nonumber\\
&  =\mathrm{div}(\tau:u)-\left(  \gamma-1\right)  \tau:\nabla u+(\gamma-1)uf.
\label{renorm_sans_eps}%
\end{align}
Next, we take the difference between $\left(  \text{\ref{renorm_sans_eps}%
}\right)  $ and $\left(  \text{\ref{passage_limite}}\right)  $ we get that%
\begin{align}
&  \partial_{t}\left(  \overline{\rho^{\gamma}}-\rho^{\gamma}\right)
+\operatorname{div}\left(  \left(  \overline{\rho^{\gamma}}-\rho^{\gamma
}\right)  u\right)  +\left(  \gamma-1\right)  \left(  \overline{\rho^{\gamma}%
}-\rho^{\gamma}\right)  \operatorname{div}u\nonumber\\
&  =-\left(  \gamma-1\right)  \left\{  \overline{\tau:\nabla u}-\tau:\nabla
u\right\}  \label{difference1}%
\end{align}
Observe that the RHS term is positive. Observe also that, formally by
multypling the above identity with $\frac{1}{\gamma}\left(  \overline
{\rho^{\gamma}}-\rho^{\gamma}\right)  ^{\frac{1}{\gamma}-1}$ the LHS of the
above expression can be written as the time-space divergence of some vector
field, see the heuristics in the introduction. The rigurous justification is a
bit more involved. First of all, the RHS of $\left(  \text{\ref{difference}%
}\right)  $ is only a measure in time and space such that we need to
regularize with respect to time and space in order to justify nonlinear change
of variables. Second of all, an even more serious problem comes from the fact
that since
\[
\partial_{t}(\rho^{\varepsilon})^{\gamma}+\operatorname{div}\left(
(\rho^{\varepsilon})^{\gamma}u^{\varepsilon}\right)  +\left(  \gamma-1\right)
(\rho^{\varepsilon})^{\gamma}\operatorname{div}u^{\varepsilon}=0
\]
and
\[
(\rho^{\varepsilon})^{\gamma}\operatorname{div}u^{\varepsilon}\text{ is
uniformly bounded in }L^{1}((0,T)\times\mathbb{T}^{3})
\]
the classical Aubin-Lions argument, see the classical argument of P.L. Lions
\cite{Li0}, Appendix $C,$ page $178$ or Lemma $6.2.$ from \cite{NoSt} allowing
to obtain that%
\[
(\rho^{\varepsilon})^{\gamma}\rightarrow\overline{\rho^{\gamma}}\text{
strongly in }C([0,T];L_{weak}^{1}(\mathbb{T}^{3}))
\]
cannot be used in this situation. To justify the formal calculation presented
in the introduction, we first prove the following

{ }

\begin{lemma}
\label{auxiliar}For any $0<s<t<T$ we have that%
\[
\frac{1}{t}\int_{0}^{t}\int_{\mathbb{T}^{3}}\left(  \overline{\rho^{\gamma}%
}\left(  \tau,x\right)  -\rho^{\gamma}\left(  \tau,x\right)  \right)
^{\frac{1}{\gamma}}d\tau dx\leq\frac{1}{s}\int_{0}^{s}\int_{\mathbb{T}^{3}%
}\left(  \overline{\rho^{\gamma}}\left(  \tau,x\right)  -\rho^{\gamma}\left(
\tau,x\right)  \right)  ^{\frac{1}{\gamma}}d\tau dx.
\]

\end{lemma}

\noindent Then in order to conclude to the identification of $\overline
{\rho^{\gamma}}$ with $\rho^{\gamma}$, we will show that
\begin{equation}
\lim_{s\rightarrow0}\frac{1}{s}\int_{0}^{s}\int_{\mathbb{T}^{3}}\left(
\overline{\rho^{\gamma}}\left(  \tau,x\right)  -\rho^{\gamma}\left(
\tau,x\right)  \right)  ^{\frac{1}{\gamma}}d\tau dx=0.
\label{initial_continuity}%
\end{equation}

\noindent\textbf{Proof of Lemma \ref{auxiliar}.} We denote by
\[
\delta\overset{not.}{=}\overline{\rho^{\gamma}}-\rho^{\gamma}\text{ }\qquad
\mu\overset{not.}{=}\overline{\tau:\nabla u}-\tau:\nabla u
\]
and thus $\left(  \text{\ref{difference1}}\right)  $ rewrites as
\begin{equation}
\partial_{t}\delta+\operatorname{div}\left(  \delta u\right)  +\left(
\gamma-1\right)  \delta\operatorname{div}u=-\left(  \gamma-1\right)
\mu\label{D_prime}%
\end{equation}
which holds true in $\mathcal{D}^{\prime}\left(  \left(  0,T\right)
\times\mathbb{T}^{3}\right)  $. Consider any $s,t\in\left(  0,T\right)  $ such
that $0<s<t<T$. Consider $n\in\mathbb{N}^{\ast}$ fixed arbitrarly such that
$\frac{1}{n}<s$. We regularize the equation $\left(  \text{\ref{D_prime}%
}\right)  $ in space-time with the help of a approximation of the identity of
the form
\[
\omega_{\varepsilon^{\prime}}\left(  t,x\right)  =\frac{1}{(\varepsilon
^{\prime})^{4}}\omega\left(  \frac{t}{\varepsilon^{\prime}}\right)
\omega\left(  \frac{\left\vert x\right\vert }{(\varepsilon^{\prime})^{3}%
}\right)  .
\]
We denote by%
\[
\delta_{\varepsilon^{\prime}}=\omega_{\varepsilon^{\prime}}\left(  t,x\right)
\ast_{t,x}\delta\text{, }\mu_{\varepsilon^{\prime}}=\omega_{\varepsilon
^{\prime}}\left(  t,x\right)  \ast_{t,x}\mu
\]
which makes sense in $\mathcal{D}^{\prime}\left(  \left(  \frac{1}%
{n},T\right)  \times\mathbb{T}^{3}\right)  $ as soon as $\varepsilon^{\prime}$
is sufficiently small. Applying $\omega_{\varepsilon^{\prime}}\left(
t,x\right)  \ast_{t,x}$ to equation $\left(  \text{\ref{D_prime}}\right)  $ we
end up with
\[
\partial_{t}\delta_{\varepsilon^{\prime}}+\operatorname{div}\left(
\delta_{\varepsilon^{\prime}}u\right)  +\left(  \gamma-1\right)
\delta_{\varepsilon^{\prime}}\operatorname{div}u=r_{\varepsilon^{\prime}%
}\left(  \delta,u\right)  -\left(  \gamma-1\right)  \mu_{\varepsilon^{\prime}}%
\]
which holds in $\mathcal{D}^{\prime}\left(  \left(  \frac{1}{n},T\right)
\times\mathbb{T}^{3}\right)  $. Above, we have that%
\begin{align}
r_{\varepsilon^{\prime}}\left(  \delta,u\right)   &  =\operatorname{div}%
\left(  (\omega_{\varepsilon^{\prime}}\ast_{t,x}\delta)u-\omega_{\varepsilon
^{\prime}}\ast_{t,x}(\delta u)\right) \label{r_eps}\\
&  +\left(  \gamma-1\right)  \left(  \left(  \omega_{\varepsilon^{\prime}}%
\ast_{t,x}\delta\right)  \operatorname{div}u-\omega_{\varepsilon^{\prime}}%
\ast_{t,x}\left(  \delta\operatorname{div}u\right)  \right)  .\nonumber
\end{align}
Since all the terms are regular, the abve equation acctually holds a.e. on
$\left(  \frac{1}{n},T\right)  \times\mathbb{T}^{3}$. We multiply the equation
with $\frac{1}{\gamma}(h+\delta_{\varepsilon^{\prime}})^{\frac{1}{\gamma}-1}$
where $h$ is a fixed positive constant. We end up with%
\begin{align*}
&  \partial_{t}\left(  h+\delta_{\varepsilon^{\prime}}\right)  ^{\frac
{1}{\gamma}}+\operatorname{div}\left(  \left(  h+\delta_{\varepsilon^{\prime}%
}\right)  ^{\frac{1}{\gamma}}u\right)  -(h+\delta_{\varepsilon^{\prime}%
})^{\frac{1}{\gamma}-1}h\operatorname{div}u\\
&  =\frac{1}{\gamma}(h+\delta_{\varepsilon^{\prime}})^{\frac{1}{\gamma}%
-1}r_{\varepsilon^{\prime}}\left(  \delta,u\right)  -\frac{1}{\gamma}%
(h+\delta_{\varepsilon^{\prime}})^{\frac{1}{\gamma}-1}\left(  \gamma-1\right)
\mu_{\varepsilon^{\prime}}.
\end{align*}
Now, consider any $\tilde{s}\in\left(  \frac{1}{n},s\right)  $ and any
$\tilde{t}\in\left(  s,t\right)  $. Let us integrate the above relation
between $\tilde{s}$ and $\tilde{t}$ in order to get that%
\begin{align*}
&  \int_{\mathbb{T}^{3}}\left(  h+\delta_{\varepsilon^{\prime}}\right)
^{\frac{1}{\gamma}}\left(  \tilde{t}\right) \\
&  =\int_{\mathbb{T}^{3}}\left(  h+\delta_{\varepsilon^{\prime}}\right)
^{\frac{1}{\gamma}}\left(  \tilde{s}\right)  +\int_{\tilde{s}}^{\tilde{t}}%
\int_{\mathbb{T}^{3}}\left[  (h+\delta_{\varepsilon^{\prime}})^{\frac
{1}{\gamma}-1}h\operatorname{div}u+\frac{1}{\gamma}(h+\delta_{\varepsilon
^{\prime}})^{\frac{1}{\gamma}-1}r_{\varepsilon^{\prime}}\left(  \delta
,u\right)  -\frac{1}{\gamma}(h+\delta_{\varepsilon^{\prime}})^{\frac{1}%
{\gamma}-1}\left(  \gamma-1\right)  \mu_{\varepsilon^{\prime}}\right] \\
&  \leq\int_{\mathbb{T}^{3}}\left(  h+\delta_{\varepsilon^{\prime}}\right)
^{\frac{1}{\gamma}}\left(  \tilde{s}\right)  +\int_{\tilde{s}}^{\tilde{t}}%
\int_{\mathbb{T}^{3}}\left[  (h+\delta_{\varepsilon^{\prime}})^{\frac
{1}{\gamma}-1}h\operatorname{div}u+\frac{1}{\gamma}(h+\delta_{\varepsilon
^{\prime}})^{\frac{1}{\gamma}-1}r_{\varepsilon^{\prime}}\left(  \delta
,u\right)  \right] \\
&  \leq\int_{\mathbb{T}^{3}}\left(  h+\delta_{\varepsilon^{\prime}}\right)
^{\frac{1}{\gamma}}\left(  \tilde{s}\right)  +\int_{\frac{1}{n}}^{T}%
\int_{\mathbb{T}^{3}}\left[  (h+\delta_{\varepsilon^{\prime}})^{\frac
{1}{\gamma}-1}h\left\vert \operatorname{div}u\right\vert +\frac{1}{\gamma
}(h+\delta_{\varepsilon^{\prime}})^{\frac{1}{\gamma}-1}\left\vert
r_{\varepsilon^{\prime}}\left(  \delta,u\right)  \right\vert \right]  .
\end{align*}
The first inequality is justified by combining the positiveness of the measure
$\mu$ (which is obtained using the lower semi-continuity assumption $\left(
\text{\ref{H2}}\right)  $) along with the fact that the convolution kernel is
positive. We integrate the above inequality with respect to $\tilde{t}$ on
$\left(  s,t\right)  $ and with respect to $\tilde{s}$ on $\left(  \frac{1}%
{n},s\right)  $ in order to recover that%
\begin{align*}
&  \left(  s-\frac{1}{n}\right)  \int_{s}^{t}\int_{\mathbb{T}^{3}}\left(
h+\delta_{\varepsilon^{\prime}}\right)  ^{\frac{1}{\gamma}}\left(  \tilde
{t}\right)  d\tilde{t}dx\leq\left(  t-s\right)  \int_{\frac{1}{n}}^{s}%
\int_{\mathbb{T}^{3}}\left(  h+\delta_{\varepsilon^{\prime}}\right)
^{\frac{1}{\gamma}}\left(  \tilde{s}\right)  d\tilde{s}dx\\
&  +\left(  t-s\right)  \left(  s-\frac{1}{n}\right)  \int_{\frac{1}{n}}%
^{T}\int_{\mathbb{T}^{3}}\left[  (h+\delta_{\varepsilon^{\prime}})^{\frac
{1}{\gamma}-1}h\left\vert \operatorname{div}u\right\vert +\frac{1}{\gamma
}(h+\delta_{\varepsilon^{\prime}})^{\frac{1}{\gamma}-1}\left\vert
r_{\varepsilon^{\prime}}\left(  \delta,u\right)  \right\vert \right]  .
\end{align*}
with $r_{\varepsilon^{\prime}}$ given by \eqref{r_eps}. We add up to the
previous inequality the quantity
\[
\left(  s-\frac{1}{n}\right)  \int_{\frac{1}{n}}^{s}\int_{\mathbb{T}^{3}%
}\left(  h+\delta_{\varepsilon^{\prime}}\right)  ^{\frac{1}{\gamma}}\left(
\tilde{s}\right)  d\tilde{s}dx
\]
which gives us%
\begin{align*}
&  \left(  s-\frac{1}{n}\right)  \int_{\frac{1}{n}}^{t}\int_{\mathbb{T}^{3}%
}\left(  h+\delta_{\varepsilon^{\prime}}\right)  ^{\frac{1}{\gamma}}\left(
\tilde{t}\right)  d\tilde{t}dx\leq\left(  t-\frac{1}{n}\right)  \int_{\frac
{1}{n}}^{s}\int_{\mathbb{T}^{3}}\left(  h+\delta_{\varepsilon^{\prime}%
}\right)  ^{\frac{1}{\gamma}}\left(  \tilde{s}\right)  d\tilde{s}dx\\
&  +\left(  t-s\right)  \left(  s-\frac{1}{n}\right)  \int_{\frac{1}{n}}%
^{T}\int_{\mathbb{T}^{3}}\left[  (h+\delta_{\varepsilon^{\prime}})^{\frac
{1}{\gamma}-1}h\left\vert \operatorname{div}u\right\vert +\frac{1}{\gamma
}(h+\delta_{\varepsilon^{\prime}})^{\frac{1}{\gamma}-1}\left\vert
r_{\varepsilon^{\prime}}\left(  \delta,u\right)  \right\vert \right]  .
\end{align*}
From the above we infer that%
\begin{align*}
\frac{1}{t-\frac{1}{n}}\int_{\frac{1}{n}}^{t}\int_{\mathbb{T}^{3}}\left(
h+\delta_{\varepsilon^{\prime}}\right)  ^{\frac{1}{\gamma}}\left(  \tilde
{t}\right)  d\tilde{t}dx  &  \leq\frac{1}{s-\frac{1}{n}}\int_{\frac{1}{n}}%
^{s}\int_{\mathbb{T}^{3}}\left(  h+\delta_{\varepsilon^{\prime}}\right)
^{\frac{1}{\gamma}}\left(  \tilde{s}\right)  d\tilde{s}dx\\
&  +\int_{\frac{1}{n}}^{T}\int_{\mathbb{T}^{3}}\left[  (h+\delta
_{\varepsilon^{\prime}})^{\frac{1}{\gamma}-1}h\left\vert \operatorname{div}%
u\right\vert +\frac{1}{\gamma}(h+\delta_{\varepsilon^{\prime}})^{\frac
{1}{\gamma}-1}\left\vert r_{\varepsilon^{\prime}}\left(  \delta,u\right)
\right\vert \right]  .
\end{align*}
Thanks to Proposition \ref{Prop_ren1}, we know that
\[
r_{\varepsilon^{\prime}}\left(  \delta,u\right)  \rightarrow0\text{ in }%
L^{1}\left(  \left(  \frac{1}{n},T\right)  \times\mathbb{T}^{3}\right)  .
\]
Observing that $(h+\delta_{\varepsilon^{\prime}})^{1/\gamma-1}\leq
h^{1/\gamma-1}$ (because $\gamma>1$ and $\delta_{\varepsilon^{\prime}}\geq0$),
we have that
\[
\int_{0}^{T}\int_{\mathbb{T}^{3}}(h+\delta_{\varepsilon^{\prime}})^{\frac
{1}{\gamma}-1}\left\vert r_{\varepsilon^{\prime}}\left(  \delta,u\right)
\right\vert \leq h^{\frac{1}{\gamma}-1}\int_{0}^{T}\int_{\mathbb{T}^{3}%
}\left\vert r_{\varepsilon^{\prime}}\left(  \delta,u\right)  \right\vert
\]
and we conclude that
\[
(h+\delta_{\varepsilon^{\prime}})^{\frac{1}{\gamma}-1}h\left\vert
\operatorname{div}u\right\vert +\frac{1}{\gamma}(h+\delta_{\varepsilon
^{\prime}})^{\frac{1}{\gamma}-1}\left\vert r_{\varepsilon^{\prime}}\left(
\delta,u\right)  \right\vert \leq\left(  1-\frac{1}{\gamma}\right)
h^{\frac{1}{\gamma}-1}\left\vert r_{\varepsilon^{\prime}}\left(
\delta,u\right)  \right\vert +h^{\frac{1}{\gamma}}\left\vert
\operatorname{div}u\right\vert .
\]
Taking into account the last observations, by making $\varepsilon^{\prime
}\rightarrow0$ we get that
\begin{align*}
&  \frac{1}{t-\frac{1}{n}}\int_{\frac{1}{n}}^{t}\int_{\mathbb{T}^{3}}\left(
h+\overline{\rho^{\gamma}}\left(  \tau,x\right)  -\rho^{\gamma}\left(
\tau,x\right)  \right)  ^{\frac{1}{\gamma}}d\tau dx\\
&  \leq\frac{1}{s-\frac{1}{n}}\int_{\frac{1}{n}}^{s}\int_{\mathbb{T}^{3}%
}\left(  h+\overline{\rho^{\gamma}}\left(  \tau,x\right)  -\rho^{\gamma
}\left(  \tau,x\right)  \right)  ^{\frac{1}{\gamma}}d\tau dx+h^{1/\gamma}%
\int_{0}^{T}\int_{{\mathbb{T}}^{3}}|\operatorname{div}\,u|.
\end{align*}
Letting $h$ go to zero we end up with
\begin{equation}
\frac{1}{t-\frac{1}{n}}\int_{\frac{1}{n}}^{t}\int_{\mathbb{T}^{3}}\left(
\overline{\rho^{\gamma}}\left(  \tau,x\right)  -\rho^{\gamma}\left(
\tau,x\right)  \right)  ^{\frac{1}{\gamma}}d\tau dx\leq\frac{1}{s-\frac{1}{n}%
}\int_{\frac{1}{n}}^{s}\int_{\mathbb{T}^{3}}\left(  \overline{\rho^{\gamma}%
}\left(  \tau,x\right)  -\rho^{\gamma}\left(  \tau,x\right)  \right)
^{\frac{1}{\gamma}}d\tau dx. \label{almost}%
\end{equation}
Since $n\in\mathbb{N}$ was chosen arbitrarly such as $\frac{1}{n}<s<t,$ we
infer that $\left(  \text{\ref{almost}}\right)  $ holds for all $n\in
\mathbb{N}$ such that $n>1/s.$ The fact that%
\[
\left(  \overline{\rho^{\gamma}}\left(  \tau,x\right)  -\rho^{\gamma}\left(
\tau,x\right)  \right)  ^{\frac{1}{\gamma}}\in L^{2\gamma}\left(  \left(
0,T\right)  \times\mathbb{T}^{3}\right)
\]
makes it possible to pass $n\rightarrow+\infty$ and to infer that%
\[
\frac{1}{t}\int_{0}^{t}\int_{\mathbb{T}^{3}}\left(  \overline{\rho^{\gamma}%
}\left(  \tau,x\right)  -\rho^{\gamma}\left(  \tau,x\right)  \right)
^{\frac{1}{\gamma}}d\tau dx\leq\frac{1}{s}\int_{0}^{s}\int_{\mathbb{T}^{3}%
}\left(  \overline{\rho^{\gamma}}\left(  \tau,x\right)  -\rho^{\gamma}\left(
\tau,x\right)  \right)  ^{\frac{1}{\gamma}}d\tau dx.
\]
This concludes the proof of Lemma \ref{auxiliar}.

\bigskip

\noindent\textbf{The final step to prove that $\overline{\rho^{\gamma}}%
=\rho^{\gamma}$.} Using Lemma \ref{auxiliar}, in order to conclude to the
identification of $\overline{\rho^{\gamma}}$ with $\rho^{\gamma}$ we only need
to show that
\begin{equation}
\lim_{s\rightarrow0}\frac{1}{s}\int_{0}^{s}\int_{\mathbb{T}^{3}}\left(
\overline{\rho^{\gamma}}\left(  \tau,x\right)  -\rho^{\gamma}\left(
\tau,x\right)  \right)  ^{\frac{1}{\gamma}}d\tau
dx=0.\label{initial_continuity_1}%
\end{equation}
\textit{In order to prove the last relation we use in a crucial manner that
the sequence of approximate solutions }$\left(  \rho^{\varepsilon}\right)
_{\varepsilon}$\textit{ verifies the energy inequality a.e in time:}
\begin{equation}
\int_{\mathbb{T}^{3}}(\rho^{\varepsilon})^{\gamma}\left(  t,x\right)
dx+\int_{0}^{t}\int_{\mathbb{T}^{3}}\tau^{\varepsilon}:\nabla u^{\varepsilon
}\leq\int_{\mathbb{T}^{3}}(\rho_{0}^{\varepsilon}(x))^{\gamma}dx+\int_{0}%
^{t}\int_{\mathbb{T}^{3}}u^{\varepsilon}f^{\varepsilon}.\label{a.e._energy}%
\end{equation}
\textit{This allows to reduce the proof of }$\left(
\text{\ref{initial_continuity_1}}\right)  $\textit{ to a continuity property
for the limit density }$\rho$\textit{. } Indeed, let us observe that $\left(
\text{\ref{a.e._energy}}\right)  $ implies that for all $s\in\left(
0,T\right)  $ we have that:%
\[
\frac{1}{s}\int_{0}^{s}\int_{\mathbb{T}^{3}}(\rho^{\varepsilon})^{\gamma
}\left(  \tau,x\right)  dxd\tau\leq\int_{\mathbb{T}^{3}}(\rho_{0}%
^{\varepsilon}\left(  x\right)  )^{\gamma}dx+\frac{1}{s}\int_{0}^{s}\left(
\int_{0}^{\tau}\int_{\mathbb{T}^{3}}u^{\varepsilon}f^{\varepsilon}\right)
d\tau.
\]
Using $\left(  \text{\ref{weak_conv1}}\right)  $ we infer that%
\[
\frac{1}{s}\int_{0}^{s}\int_{\mathbb{T}^{3}}\overline{\rho^{\gamma}}\left(
\tau,x\right)  dxd\tau\leq\int_{\mathbb{T}^{3}}\rho_{0}^{\gamma}(x)dx+\frac
{1}{s}\int_{0}^{s}\left(  \int_{0}^{\tau}\int_{\mathbb{T}^{3}}uf\right)
d\tau.
\]
Next, we use H\"{o}lder's inequality to infer that%
\begin{align*}
&  \frac{1}{s}\int_{0}^{s}\int_{\mathbb{T}^{3}}\left(  \overline{\rho^{\gamma
}}\left(  \tau,x\right)  -\rho^{\gamma}\left(  \tau,x\right)  \right)
^{\frac{1}{\gamma}}d\tau dx\\
&  \leq\frac{1}{s}\int_{0}^{s}\int_{\mathbb{T}^{3}}\left(  \overline
{\rho^{\gamma}}\left(  \tau,x\right)  -\rho^{\gamma}\left(  \tau,x\right)
\right)  d\tau dx\\
&  =\frac{1}{s}\int_{0}^{s}\int_{\mathbb{T}^{3}}\left(  \overline{\rho
^{\gamma}}\left(  \tau,x\right)  -\rho_{0}^{\gamma}\left(  x\right)  \right)
d\tau dx+\frac{1}{s}\int_{0}^{s}\int_{\mathbb{T}^{3}}\left(  \rho_{0}^{\gamma
}\left(  x\right)  -\rho^{\gamma}\left(  \tau,x\right)  \right)  d\tau dx\\
&  \leq\frac{1}{s}\int_{0}^{s}\left(  \int_{0}^{\tau}\int_{\mathbb{T}^{3}%
}uf\right)  d\tau+\frac{1}{s}\int_{0}^{s}\int_{\mathbb{T}^{3}}\left(  \rho
_{0}^{\gamma}\left(  x\right)  -\rho^{\gamma}\left(  \tau,x\right)  \right)
d\tau dx\\
&  =\frac{1}{s}\int_{0}^{s}\left(  \int_{0}^{\tau}\int_{\mathbb{T}^{3}%
}uf\right)  d\tau+\int_{\mathbb{T}^{3}}\rho_{0}^{\gamma}\left(  x\right)
dx-\frac{1}{s}\int_{0}^{s}\int_{\mathbb{T}^{3}}\rho^{\gamma}\left(
\tau,x\right)  d\tau dx.
\end{align*}
Thus, since
\[
uf\in L^{1}\left(  \left(  0,T\right)  \times\mathbb{T}^{3}\right)
\]
proving $\left(  \text{\ref{initial_continuity}}\right)  $ reduces to prove
that
\[
\lim_{s\rightarrow0}\frac{1}{s}\int_{0}^{s}\int_{\mathbb{T}^{3}}\left(
\rho_{0}^{\gamma}\left(  x\right)  -\rho^{\gamma}\left(  \tau,x\right)
\right)  d\tau dx=0.
\]
The proof of the above is contained in the following

\begin{lemma}
\label{transport_lemma}Consider $\rho\in L^{\infty}\left(  \left(  0,T\right)
;L^{\gamma}\left(  \mathbb{T}^{3}\right)  \right)  \cap C([0,T];L_{weak}%
^{\gamma}(\mathbb{T}^{3}))\cap L^{2\gamma}\left(  \left(  0,T\right)
\times\mathbb{T}^{3}\right)  $ and $u\in L^{2}\left(  \left(  0,T\right)
;H^{1}\left(  \mathbb{T}^{3}\right)  \right)  $ verifying the transport
equation
\[
\partial_{t}\rho+\operatorname{div}\left(  \rho u\right)  =0\text{ in
}\mathcal{D}^{\prime}\left(  \left(  0,T\right)  \times\mathbb{T}^{3}\right)
\]
along with the fact that%
\[
\lim_{t\rightarrow0}\int_{\mathbb{T}^{3}}\rho\left(  t,x\right)  \psi\left(
x\right)  dx=\int_{\mathbb{T}^{3}}\rho_{0}\left(  x\right)  \psi\left(
x\right)  dx\text{ for all }\psi\in C_{per}^{\infty}\left(  \mathbb{R}%
^{d}\right)  .
\]
Then
\begin{equation}
\lim_{s\rightarrow0}\frac{1}{s}\int_{0}^{s}\int_{\mathbb{T}^{3}}\left(
\rho^{\gamma}\left(  \tau,x\right)  -\rho_{0}^{\gamma}\left(  x\right)
\right)  d\tau dx=0. \label{time_mean_continuity}%
\end{equation}
\bigskip
\end{lemma}

\noindent\textbf{Proof of Lemma \ref{transport_lemma}.} First of all, it is
classical to recover that $\rho\in C\left(  [0,T);L^{p}\left(  \mathbb{T}%
^{3}\right)  \right)  $ with $p\in\lbrack1,\gamma)$ and that%
\begin{equation}
\lim_{t\rightarrow0}\rho\left(  t,\right)  =\rho_{0}\text{ in }L^{p}\text{ for
all }p\in\lbrack1,\gamma). \label{weak_time_cont}%
\end{equation}
This is of course not sufficient in order to prove $\left(
\text{\ref{time_mean_continuity}}\right)  $. Let us consider a spatial
approximation of the identity $\left(  \omega_{\varepsilon}\right)
_{\varepsilon>0}=\left(  \frac{1}{\varepsilon^{3}}\omega\left(  \frac{\cdot
}{\varepsilon}\right)  \right)  _{\varepsilon>0}$. We will denote by
\[
\rho_{\varepsilon}\left(  t,x\right)  =\omega_{\varepsilon}\ast\rho\left(
t,x\right)  .
\]
We have that%
\[
\lim_{\varepsilon\rightarrow0}\left\Vert \rho-\rho_{\varepsilon}\right\Vert
_{L^{2\gamma}((0,T)\times\mathbb{T}^{3})}=0.
\]
Moreover, using \ref{weak_time_cont} for all $\varepsilon>0$ we have that%
\begin{equation}
\lim_{t\rightarrow0}\rho_{\varepsilon}\left(  t,\cdot\right)  =\omega
_{\varepsilon}\ast\rho_{0}\text{ in }L^{\gamma}. \label{weak_time_continuity}%
\end{equation}
For example,
\[
\left\Vert \rho_{\varepsilon}\left(  t,\cdot\right)  -\omega_{\varepsilon}%
\ast\rho_{0}\right\Vert _{L^{\gamma}}\leq\left\Vert \omega_{\varepsilon
}\right\Vert _{L^{p\left(  \eta\right)  }(\mathbb{T}^{3})}\left\Vert
\rho\left(  t,\cdot\right)  -\rho_{0}\right\Vert _{L^{\gamma-\eta}%
(\mathbb{T}^{3})}.
\]
Next, we apply $\omega_{\varepsilon}$ for the transport equation such as to
obtain
\begin{equation}
\partial_{t}\rho_{\varepsilon}^{\gamma}+\operatorname{div}\left(
\rho_{\varepsilon}^{\gamma}u\right)  +\left(  \gamma-1\right)  \rho
_{\varepsilon}^{\gamma}\operatorname{div}u=\gamma\rho_{\varepsilon}^{\gamma
-1}r_{\varepsilon}\text{ in }\mathcal{D}^{\prime}\left(  \left(  0,T\right)
\times\mathbb{T}^{3}\right)  \label{asterix}%
\end{equation}
with
\[
r_{\varepsilon}\rightarrow0\text{ in }L^{\frac{2\gamma}{\gamma+1}}\left(
\left(  0,T\right)  \times\mathbb{T}^{3}\right)  .
\]
An important property is that for all $\varepsilon>0$ and a.e. $t\in\left(
0,T\right)  $ it holds true that
\begin{align}
h_{\varepsilon}\left(  t\right)   &  =\int_{\mathbb{T}^{3}}\gamma
\rho_{\varepsilon}^{\gamma-1}\left(  t\right)  r_{\varepsilon}\left(
t\right)  -\left(  \gamma-1\right)  \rho_{\varepsilon}^{\gamma}\left(
t\right)  \operatorname{div}u\left(  t\right)  .\nonumber\\
&  \leq\left(  \gamma-1\right)  \int_{\mathbb{T}^{3}}\rho_{\varepsilon
}^{\gamma}\left(  t\right)  \left\vert \operatorname{div}u\left(  t\right)
\right\vert +\gamma\int_{\mathbb{T}^{3}}\rho_{\varepsilon}^{\gamma
-1}\left\vert r_{\varepsilon}\right\vert \nonumber\\
&  \leq C_{\gamma}\left\Vert \rho\left(  t\right)  \right\Vert _{L^{2\gamma
}(\mathbb{T}^{3})}^{\gamma}\left\Vert \nabla u\left(  t\right)  \right\Vert
_{L^{2}(\mathbb{T}^{3})}:=h\left(  t\right)  \in L^{1}\left(  0,T\right)  .
\label{important_estimate}%
\end{align}
Integrating the $\left(  \text{\ref{asterix}}\right)  $ we end up with%
\[
\frac{d}{dt}\int_{\mathbb{T}^{3}}\rho_{\varepsilon}^{\gamma}\left(
t,x\right)  dx=h_{\varepsilon}\left(  t\right)  \in L^{1}\left(  0,T\right)
.
\]
But using $\left(  \text{\ref{weak_time_continuity}}\right)  $ along with the
last relation we obtain that the application $t\rightarrow$ $\int
_{\mathbb{T}^{3}}\rho_{\varepsilon}^{\gamma}\left(  t\right)  $ is absolutely
continious and we may write that%
\[
\int_{\mathbb{T}^{3}}\rho_{\varepsilon}^{\gamma}\left(  t,x\right)
dx=\int_{\mathbb{T}^{3}}(\omega_{\varepsilon}\ast\rho_{0})^{\gamma}\left(
x\right)  dx+\int_{0}^{t}h_{\varepsilon}\left(  \tau\right)  d\tau.
\]
From this and $\left(  \text{\ref{important_estimate}}\right)  $ we learn that%
\[
\left\vert \int_{\mathbb{T}^{3}}\rho_{\varepsilon}^{\gamma}\left(  t,x\right)
dx-\int_{\mathbb{T}^{3}}(\omega_{\varepsilon}\ast\rho_{0})^{\gamma}\left(
x\right)  dx\right\vert \leq\int_{0}^{t}h\left(  \tau\right)  d\tau.
\]
Now, we know that $h\left(  t\right)  \in L^{1}\left(  0,T\right)  $ and
consequently the application $t\rightarrow\int_{0}^{t}h\left(  \tau\right)
d\tau$ is absolutely continious and
\[
\lim_{t\rightarrow0}\int_{0}^{t}h\left(  \tau\right)  d\tau=0.
\]
Let us fix $\eta>0$. Using the above we obtain the existence of a $t_{\eta}>0$
such that for all $t\in\left(  0,t_{\eta}\right)  $ \textit{and for all
}$\varepsilon>0$ one has%
\[
\left\vert \int_{\mathbb{T}^{3}}\rho_{\varepsilon}^{\gamma}\left(  t,x\right)
dx-\int_{\mathbb{T}^{3}}(\omega_{\varepsilon}\ast\rho_{0})^{\gamma}\left(
x\right)  dx\right\vert \leq\int_{0}^{t}h\left(  \tau\right)  d\tau\leq\eta.
\]
By the triangle inequality, we have that for all $\varepsilon>0$ and
$t\in\left(  0,t_{\eta}\right)  $%
\begin{equation}
\left\vert \frac{1}{t}\int_{0}^{t}\int_{\mathbb{T}^{3}}\rho_{\varepsilon
}^{\gamma}\left(  \tau,x\right)  dxd\tau-\int_{\mathbb{T}^{3}}(\omega
_{\varepsilon}\ast\rho_{0})^{\gamma}\left(  x\right)  dx\right\vert \leq\eta.
\label{almost_almost}%
\end{equation}
For $t$ fixed arbitrarly in $\left(  0,t_{\eta}\right)  $ we use the fact that%
\[
\lim_{\varepsilon\rightarrow0}\left\Vert \rho-\rho_{\varepsilon}\right\Vert
_{L^{2\gamma}((0,T)\times\mathbb{T}^{3})}=0
\]
we pass to the limit into $\left(  \text{\ref{almost_almost}}\right)  $ in
order to obtain that for all $t\in\left(  0,t_{\eta}\right)  $%
\[
\left\vert \frac{1}{t}\int_{0}^{t}\int_{\mathbb{T}^{3}}\rho^{\gamma}\left(
\tau,x\right)  dxd\tau-\int_{\mathbb{T}^{3}}\rho_{0}^{\gamma}\left(  x\right)
dx\right\vert \leq\eta.
\]
Since $\eta$ was fixed arbitrarly, the last property translates that
\[
\lim_{t\rightarrow0}\frac{1}{t}\int_{0}^{t}\int_{\mathbb{T}^{3}}\rho^{\gamma
}\left(  \tau,x\right)  dxd\tau=\int_{\mathbb{T}^{3}}\rho_{0}^{\gamma}\left(
x\right)  dx.
\]
This concludes the proof of Lemma \ref{transport_lemma}.

\noindent Using Lemma \ref{auxiliar} and the limit property
\eqref{initial_continuity}, we conclude that
\[
\overline{\rho^{\gamma}}=\rho^{\gamma}\text{ a.e. on }\left(  0,T\right)
\times\mathbb{T}^{3}.
\]

\section{Construction of solutions\label{Construction}}

In this section, we propose a regularized system with diffusion and drag terms
on the density for which we prove global existence and uniqueness of strong
solution on $(0,T)$ using a fixed point procedure. Then passing to the limit
with respect to the regularization parameter provides a global solution of the
quasi-stationary compressible Stokes system with diffusion on the density and
drag terms on the density. It remains to show that these extra terms do not
perturb the stability procedure, we explained in subsection \ref{weakstab}, to
prove Theorem \ref{main}.

\subsection{The approximate system}

Let us be more precise. For any fixed strictly positive parameter
$\varepsilon,\delta$ we are able to construct a global solution of the
following regularized version of the original system:%
\begin{equation}
\left\{
\begin{array}
[c]{l}%
\partial_{t}\rho+\operatorname{div}\left(  \rho\omega_{\delta}\ast u\right)
=\varepsilon\Delta\rho-\varepsilon\rho^{2\gamma}-\varepsilon\rho^{2\gamma
+1}-\varepsilon\rho^{3},\\
\mathcal{A}u+\nabla\omega_{\delta}\ast\rho^{\gamma}=f,\\
\rho_{|t=0}=\rho_{0}^{reg},
\end{array}
\right.  \tag{$\mathcal{S}_\varepsilon,\delta$}\label{approx}%
\end{equation}
with $\omega_{\delta}$ the standard regularizing kernel see $\left(
\text{\ref{notation_approx}}\right)  $. The function $\rho_{0}^{reg}$ is
supposed to be regular enough as to ensure existence of solutions to the
transport equation with regular velocity and initial data initial data
$\rho_{0}^{reg}$. The construction of solutions for $\left(
\text{\ref{approx}}\right)  $ is achieved by a classical fixed point argument.

In a second time, we show that a sequence of solutions $\left(  \rho
^{\varepsilon,\delta},u^{\varepsilon,\delta}\right)  $ of $\left(
\text{\ref{approx}}\right)  $ tends, when we let $\delta$ go to zero, to
$\left(  \rho^{\varepsilon},u^{\varepsilon}\right)  $ which is a solution of
the system
\begin{equation}
\left\{
\begin{array}
[c]{l}%
\partial_{t}\rho+\operatorname{div}\left(  \rho u\right)  =\varepsilon
\Delta\rho-\varepsilon\rho^{2\gamma}-\varepsilon\rho^{2\gamma+1}%
-\varepsilon\rho^{3},\\
\mathcal{A}u+\nabla\rho^{\gamma}=f,\\
\rho_{|t=0}=\rho_{0},
\end{array}
\right.  \tag{$\mathcal{S}_\varepsilon$}\label{approx_1level}%
\end{equation}
which, moreover, verifies the following estimates, uniformly in $\varepsilon$
(we skip the $\varepsilon$ script in the inequalities bellow such to render
them more readable):%
\begin{equation}
\left\{
\begin{array}
[c]{l}%
{\displaystyle\int_{\mathbb{T}^{3}}}
\rho\left(  t\right)  +\varepsilon%
{\displaystyle\int_{0}^{t}}
{\displaystyle\int_{\mathbb{T}^{3}}}
\rho^{2\gamma}+\varepsilon%
{\displaystyle\int_{0}^{t}}
{\displaystyle\int_{\mathbb{T}^{3}}}
\rho^{2\gamma+1}+\varepsilon%
{\displaystyle\int_{0}^{t}}
{\displaystyle\int_{\mathbb{T}^{3}}}
\rho^{3}=%
{\displaystyle\int_{\mathbb{T}^{3}}}
\rho_{0},\\%
{\displaystyle\int_{\mathbb{T}^{3}}}
\rho^{\gamma}\left(  t\right)  +\frac{c\left(  \gamma-1\right)  }{2}%
{\displaystyle\int_{0}^{t}}
{\displaystyle\int_{\mathbb{T}^{3}}}
\left\vert \nabla u\right\vert ^{2}\\
\hskip2cm+\varepsilon\gamma%
{\displaystyle\int_{0}^{t}}
{\displaystyle\int_{\mathbb{T}^{3}}}
\rho^{3\gamma-1}+\varepsilon\,\gamma%
{\displaystyle\int_{0}^{t}}
{\displaystyle\int_{\mathbb{T}^{3}}}
\rho^{3\gamma}+\varepsilon\gamma%
{\displaystyle\int_{0}^{t}}
{\displaystyle\int_{\mathbb{T}^{3}}}
\rho^{\gamma+2}\\
\hskip4cm+4\varepsilon\lbrack1-\frac{1}{\gamma}]%
{\displaystyle\int_{0}^{t}}
{\displaystyle\int_{\mathbb{T}^{3}}}
\left\vert \nabla\rho^{\frac{\gamma}{2}}\right\vert ^{2}\leq C\left(
c,\gamma\right)  \left(
{\displaystyle\int_{\mathbb{T}^{3}}}
\rho_{0}^{\gamma}+\left\Vert f\right\Vert _{L_{t}^{2}L^{\frac{6}{5}}}%
^{2}\right)  ,\\
\left\Vert \rho^{\gamma}\right\Vert _{L^{2}((0,T)\times{\mathbb{T}^{3}})}\leq
C\left(  c,\gamma\right)  \left(  \sqrt{t}+\max\left\{  1,\left\Vert
A\right\Vert _{L^{\infty}}\right\}  \right)  \left(  \left\Vert \rho
_{0}\right\Vert _{L^{\gamma}}^{\frac{\gamma}{2}}+\left\Vert f\right\Vert
_{L_{t}^{2}L^{\frac{6}{5}}}\right)  ,
\end{array}
\right.
\end{equation}
with $c$ defined by $\left(  \text{\ref{H4}}\right)  $ and $C\left(
c,\gamma\right)  $ a constant depending only on $c$ and $\gamma$.

Finally, we show that we can adapt the proof of Theorem \ref{main} in order to
pass to the limit $\varepsilon\rightarrow0$ and thus obtaining a solution for
the compressible Stokes system.

\subsection{Construction of solutions for the regularized system $\left(
\text{\ref{approx}}\right)  $}

We consider $T>0$ to be precised later and we denote by
\[
L^{2}(0,T;\dot{H}^{1}(\mathbb{T}^{3}))=\left\{  u\in L^{2}(0,T;H^{1}%
(\mathbb{T}^{3})):\int_{\mathbb{T}^{3}}u\left(  t\right)  =0\text{ a.e. }%
t\in(0,T)\right\}
\]
Consider
\[
B:L^{2}(0,T;\dot{H}^{1}(\mathbb{T}^{3}))\rightarrow L^{2}(0,T;\dot{H}%
^{1}(\mathbb{T}^{3}))
\]
defined as%

\begin{equation}
\left\{
\begin{array}
[c]{l}%
\partial_{t}\rho+\operatorname{div}\left(  \rho\omega_{\delta}\ast v\right)
=\varepsilon\Delta\rho-\varepsilon\rho^{2\gamma}-\varepsilon\rho^{2\gamma
+1}-\varepsilon\rho^{3},\\
\mathcal{A}B(v)+\nabla\omega_{\delta}\ast\rho^{\gamma}=f,\\
\rho_{|t=0}=\rho_{0}^{reg}%
\end{array}
\right.  \label{point_fixe}%
\end{equation}
Obviously if $v\in L^{2}(0,T;\dot{H}^{1}(\mathbb{T}^{3}))$ then $\omega
_{\delta}\ast v\in L^{2}(0,T;C^{\infty}(\mathbb{T}^{3}))$ such that the
existence of a regular\emph{ positive} solution for the first equation of
system $\left(  \text{\ref{point_fixe}}\right)  $ follows by classical
arguments. Also, $B\left(  v\right)  $ is well-defined as an element of
$L^{2}(0,T;\dot{H}^{1}(\mathbb{T}^{3}))$ and
\[
\int_{0}^{T}\int_{\mathbb{T}^{3}}A(t,x)D(B(v):D(B(v))=\int_{0}^{T}%
\int_{\mathbb{T}^{3}}\omega_{\delta}\ast\rho^{\gamma}%
\mathrm{\operatorname{div}}B(v)+\int_{0}^{T}\int_{\mathbb{T}^{3}}fu
\]
which provides
\begin{equation}
\left\Vert \nabla B\left(  v\right)  \right\Vert _{L^{2}((0,T)\times
{\mathbb{T}}^{3})}\leq C\left\Vert \omega_{\delta}\ast\rho^{\gamma}\right\Vert
_{L^{2}((0,T)\times{\mathbb{T}}^{3})}+C\left\Vert f\right\Vert _{L_{t}%
^{2}L^{\frac{6}{5}}}, \label{B(v)_bound_1}%
\end{equation}
with $C$ depending only on the dissipation operator. Let us integrate the
equation defining $\rho$ in order to see that
\[
\int_{{\mathbb{T}}^{3}}\rho\left(  t\right)  +\varepsilon\int_{0}^{t}%
\int_{{\mathbb{T}}^{3}}\rho^{2\gamma}+\varepsilon\int_{0}^{t}\int
_{{\mathbb{T}}^{3}}\rho^{2\gamma+1}+\varepsilon\int_{0}^{t}\int_{{\mathbb{T}%
}^{3}}\rho^{3}=\int_{{\mathbb{T}}^{3}}\rho_{0}^{reg}%
\]
which, enables us to conclude, that%
\begin{equation}
\left\Vert \nabla B\left(  v\right)  \right\Vert _{L^{2}((0,T)\times
{\mathbb{T}}^{3})}\leq\tilde{C}\left(  c,\gamma\right)  \left(  \frac
{1}{\varepsilon}%
{\displaystyle\int_{\mathbb{T}^{3}}}
\rho_{0}^{reg}+\left\Vert f\right\Vert _{L_{t}^{2}L^{\frac{6}{5}}}^{2}\right)
^{\frac{1}{2}}. \label{operator_Boundness}%
\end{equation}
Thus, we conclude that for any $T>0$, the operator $B$ (trivially) maps
$E_{T}$ into itself where
\[
E_{T}=\left\{  v\in L_{T}^{2}(\dot{H}^{1}(\mathbb{T}^{3})):\left\Vert \nabla
v\right\Vert _{L^{2}((0,T)\times{\mathbb{T}}^{3})}\leq\tilde{C}\left(
c,\gamma\right)  \left(  \frac{1}{\varepsilon}%
{\displaystyle\int_{\mathbb{T}^{3}}}
\rho_{0}^{reg}+\left\Vert f\right\Vert _{L_{t}^{2}L^{\frac{6}{5}}}^{2}\right)
\right\}  .
\]
In the following, we aim at showing that $B$ is a contraction on $E_{T}$.

The first observation that we make in towards this direction is that using a
maximum principle we get%
\begin{align}
\left\Vert \rho\right\Vert _{L^{\infty}((0,t)\times{{\mathbb{T}}^{3}})}  &
\leq\left\Vert \rho_{0}^{reg}\right\Vert _{L^{\infty}({{\mathbb{T}}^{3}})}%
\exp\left(  \int_{0}^{t}\left\Vert \operatorname{div}\omega_{\delta}\ast
v\right\Vert _{L^{\infty}({{\mathbb{T}}^{3}})}\right) \nonumber\\
&  \leq\left\Vert \rho_{0}^{reg}\right\Vert _{L^{\infty}({{\mathbb{T}}^{3}}%
)}\exp\left(  \sqrt{t}C_{\varepsilon,\delta}\right)  . \label{rho_L_infty}%
\end{align}
Next, let us multiply the first equation of $\left(  \text{\ref{point_fixe}%
}\right)  $ with $\rho$ and integrate in order to obtain that%
\[
\frac{1}{2}\int_{{\mathbb{T}}^{3}}\rho^{2}+\varepsilon\int_{0}^{t}%
\int_{{\mathbb{T}}^{3}}\left\vert \nabla\rho\right\vert ^{2}+\varepsilon
\int_{0}^{t}\int_{{\mathbb{T}}^{3}}\rho^{2\gamma+1}+\varepsilon\int_{0}%
^{t}\int_{{\mathbb{T}}^{3}}\rho^{2\gamma+2}+\varepsilon\int_{0}^{t}%
\int_{{\mathbb{T}}^{3}}\rho^{4}=\varepsilon\int_{{\mathbb{T}}^{3}}\rho
^{2}\operatorname{div}\left(  \omega_{\delta}\ast v\right)
\]
and thus by Gronwall's lemma we get that%
\begin{align}
\frac{1}{2}\int_{{\mathbb{T}}^{3}}\rho^{2}+\varepsilon\int_{0}^{t}%
\int_{{\mathbb{T}}^{3}}\left\vert \nabla\rho\right\vert ^{2}+\varepsilon
\int_{0}^{t}\int_{{\mathbb{T}}^{3}}\rho^{2\gamma+1}  &  +\varepsilon\int
_{0}^{t}\int_{{\mathbb{T}}^{3}}\rho^{2\gamma+2}+\varepsilon\int_{0}^{t}%
\int_{{\mathbb{T}}^{3}}\rho^{4}\nonumber\\
&  \leq\frac{1}{2}\int_{{\mathbb{T}}^{3}}\left(  \rho_{0}^{reg}\right)
^{2}\exp\left(  \int_{0}^{t}\left\Vert \operatorname{div}\left(
\omega_{\delta}\ast v\right)  \right\Vert _{L^{\infty}({{\mathbb{T}}^{3}}%
)}\right) \nonumber\\
&  \leq\frac{1}{2}\int_{{\mathbb{T}}^{3}}\left(  \rho_{0}^{reg}\right)
^{2}\exp\left(  tC_{\delta}\int_{0}^{t}\left\Vert \nabla v\right\Vert
_{L^{2}({{\mathbb{T}}^{3}})}^{2}\right) \nonumber\\
&  \leq\frac{1}{2}\int_{{\mathbb{T}}^{3}}\left(  \rho_{0}^{reg}\right)
^{2}\exp\left(  tC_{\delta}\left(  \frac{1}{\varepsilon}%
{\displaystyle\int_{\mathbb{T}^{3}}}
\rho_{0}^{reg}+\left\Vert f\right\Vert _{L_{t}^{2}L^{\frac{6}{5}}}^{2}\right)
^{\frac{1}{2}}\right)  \label{contraction_0}%
\end{align}
Let us consider $v_{1},v_{2}\in E_{T}$ and let us consider%
\[
\left\{
\begin{array}
[c]{l}%
\partial_{t}\rho_{i}+\operatorname{div}\left(  \rho_{i}\omega_{\delta}\ast
v_{i}\right)  =\varepsilon\Delta\rho_{i}-\varepsilon\rho_{i}^{2\gamma
}-\varepsilon\rho_{i}^{2\gamma+1}-\varepsilon\rho_{i}^{3},\\
\mathcal{A}B(v_{i})+\nabla\omega_{\delta}\ast\rho_{i}^{\gamma}=0,\\
\rho_{i|t=0}=\rho_{0}^{reg}%
\end{array}
\right.
\]
with $i\in1,2$. Of course, $\rho_{1}$ and $\rho_{2}$ verify the estimate
$\left(  \text{\ref{contraction_0}}\right)  $. We denote by $r=\rho_{1}%
-\rho_{2}$ and $w=v_{1}-v_{2}$. We infer that%
\[
\left\{
\begin{array}
[c]{l}%
\partial_{t}r+\operatorname{div}\left(  r\omega_{\delta}\ast v_{1}\right)
=\varepsilon\Delta r-\varepsilon\left(  \rho_{1}^{2\gamma}+\rho_{1}%
^{2\gamma+1}+\rho_{1}^{3}-\rho_{2}^{2\gamma}-\rho_{2}^{2\gamma+1}-\rho_{2}%
^{3}\right)  -\operatorname{div}\left(  \rho_{2}V_{\delta}\ast w\right)  ,\\
\mathcal{A}\left(  B(v_{1})-B\left(  v_{2}\right)  \right)  +\nabla
\omega_{\delta}\ast\left(  \rho_{1}^{\gamma}-\rho_{2}^{\gamma}\right)  =0,\\
r_{|t=0}=0
\end{array}
\right.
\]
By multiplying the first equation with $r$ we get that%
\begin{align}
&  \int_{{\mathbb{T}}^{3}}\frac{r^{2}\left(  t\right)  }{2}+\varepsilon
\int_{0}^{t}\int_{{\mathbb{T}}^{3}}\left\vert \nabla r\right\vert
^{2}+\varepsilon\int_{0}^{t}\int_{{\mathbb{T}}^{3}}\left(  \rho_{1}^{2\gamma
}+\rho_{1}^{2\gamma+1}+\rho_{1}^{3}-\rho_{2}^{2\gamma}-\rho_{2}^{2\gamma
+1}-\rho_{2}^{3}\right)  r\nonumber\\
&  \leq\int_{0}^{t}\int_{{\mathbb{T}}^{3}}r^{2}\operatorname{div}%
\omega_{\delta}\ast v_{1}+\int_{0}^{t}\int_{{\mathbb{T}}^{3}}%
\operatorname{div}\left(  \rho_{2}\omega_{\delta}\ast w\right)  r\nonumber\\
&  \leq\int_{0}^{t}\int r^{2}\left\Vert \operatorname{div}\omega_{\delta}\ast
v_{1}\right\Vert _{L^{\infty}({{\mathbb{T}}^{3}})}+\frac{1}{2\varepsilon}%
\int_{0}^{t}\left\Vert \rho_{2}\right\Vert _{L^{2}({{\mathbb{T}}^{3}})}%
^{2}\left\Vert \omega_{\delta}\ast\delta v\right\Vert _{L^{\infty
}({{\mathbb{T}}^{3}})}^{2}+\frac{\varepsilon}{2}\int_{0}^{t}\int_{{\mathbb{T}%
}^{3}}\left\vert \nabla r\right\vert ^{2}\nonumber\\
&  \leq\int_{0}^{t}\int_{{\mathbb{T}}^{3}}r^{2}\left\Vert \operatorname{div}%
\omega_{\delta}\ast v_{1}\right\Vert _{L^{\infty}({{\mathbb{T}}^{3}}%
)}+C_{\delta,\varepsilon}\exp\left(  tC_{\delta,\varepsilon}\int\rho_{0}%
^{reg}\right)  \int_{0}^{t}\left\Vert \delta v\right\Vert _{L^{6}%
({{\mathbb{T}}^{3}})}^{2}+\frac{\varepsilon}{2}\int_{0}^{t}\int_{{\mathbb{T}%
}^{3}}\left\vert \nabla r\right\vert ^{2}\nonumber\\
&  \leq\int_{0}^{t}\int_{{\mathbb{T}}^{3}}r^{2}\left\Vert \operatorname{div}%
\omega_{\delta}\ast v_{1}\right\Vert _{L^{\infty}({{\mathbb{T}}^{3}}%
)}+C_{\delta,\varepsilon}\exp\left(  tC_{\delta,\varepsilon}\int\rho_{0}%
^{reg}\right)  \int_{0}^{t}\left\Vert \nabla\delta v\right\Vert _{L^{2}%
({{\mathbb{T}}^{3}})}^{2}+\frac{\varepsilon}{2}\int_{0}^{t}\int_{{\mathbb{T}%
}^{3}}\left\vert \nabla r\right\vert ^{2}%
\end{align}
and thus using Gr\"{o}nwall's lemma we get that%
\begin{align}
&  \int_{{\mathbb{T}}^{3}}\frac{r^{2}\left(  t\right)  }{2}+\frac{\varepsilon
}{2}\int_{0}^{t}\int_{{\mathbb{T}}^{3}}\left\vert \nabla r\right\vert
^{2}+\varepsilon\int_{0}^{t}\int_{{\mathbb{T}}^{3}}\left(  \rho_{1}^{2\gamma
}-\rho_{2}^{2\gamma}\right)  r+\varepsilon\int_{0}^{t}\int_{{\mathbb{T}}^{3}%
}\left(  \rho_{1}^{2\gamma+1}-\rho_{2}^{2\gamma+1}\right)  r+\varepsilon
\int_{0}^{t}\int_{{\mathbb{T}}^{3}}\left(  \rho_{1}^{3}-\rho_{2}^{3}\right)
r\nonumber\\
&  \leq C_{\delta,\varepsilon}\exp\left(  tC_{\delta,\varepsilon}\left(
{\displaystyle\int_{\mathbb{T}^{3}}}
\rho_{0}^{reg}+\left\Vert f\right\Vert _{L_{t}^{2}L^{\frac{6}{5}}}^{2}\right)
^{\frac{1}{2}}\right)  \int_{0}^{t}\left\Vert \nabla w\right\Vert
_{L^{2}({{\mathbb{T}}^{3}})}^{2}\exp\left(  \int_{0}^{t}\int_{{\mathbb{T}}%
^{3}}\left\Vert \operatorname{div}\omega_{\delta}\ast v_{1}\right\Vert
_{L^{\infty}({{\mathbb{T}}^{3}})}\right) \nonumber\\
&  \leq C_{\delta,\varepsilon}\exp\left(  C_{\delta,\varepsilon}t\right)
\int_{0}^{t}\left\Vert \nabla w\right\Vert _{L^{2}({{\mathbb{T}}^{3}})}%
^{2}=C_{\delta,\varepsilon}\exp\left(  C_{\delta,\varepsilon}t\right)
\int_{0}^{t}\left\Vert \nabla v_{1}-\nabla v_{2}\right\Vert _{L^{2}%
({{\mathbb{T}}^{3}})}^{2} \label{contraction_4}%
\end{align}
Finally, recalling that%
\[
\mathcal{A}\left(  B(v_{1})-B\left(  v_{2}\right)  \right)  +\nabla
\omega_{\delta}\ast\left(  \rho_{1}^{\gamma}-\rho_{2}^{\gamma}\right)  =0,
\]
we infer that%
\begin{equation}
\left\Vert \nabla\left(  B\left(  v_{1}\right)  -B(v_{2})\right)  \right\Vert
_{L^{2}((0,t)\times{{\mathbb{T}}^{3}})}\leq Ct^{\frac{1}{2}}\left\Vert
\rho_{1}^{\gamma}-\rho_{2}^{\gamma}\right\Vert _{L^{\infty}(0,t;L^{2}%
({{\mathbb{T}}^{3}}))} \label{contraction_5}%
\end{equation}
We use the intermediate value theorem and estimate $\left(
\text{\ref{rho_L_infty}}\right)  $ in order to asses that
\begin{align}
\left\vert \rho_{1}^{\gamma}-\rho_{2}^{\gamma}\right\vert  &  \leq
\gamma\left\vert \rho_{1}-\rho_{2}\right\vert \max\left\{  \left\Vert \rho
_{1}\right\Vert _{L^{\infty}((0,t)\times{{\mathbb{T}}^{3}})}^{\gamma
-1},\left\Vert \rho_{2}\right\Vert _{L^{\infty}((0,t)\times{{\mathbb{T}}^{3}%
})}^{\gamma-1}\right\} \nonumber\\
&  \leq\gamma\left\vert \rho_{1}-\rho_{2}\right\vert \left\Vert \rho_{0}%
^{reg}\right\Vert _{L^{\infty}({{\mathbb{T}}^{3}})}^{\gamma-1}\exp\left(
\sqrt{t}C_{\delta,\varepsilon}\right)
\end{align}
which, in turn implies that%
\[
\left\Vert \rho_{1}^{\gamma}-\rho_{2}^{\gamma}\right\Vert _{L^{\infty
}(0,t;L^{2}({\mathbb{T}}^{3}))}\leq\gamma\left\Vert \rho_{0}^{reg}\right\Vert
_{L^{\infty}({\mathbb{T}}^{3}))}^{\gamma-1}\exp\left(  \sqrt{t}C_{\delta
,\varepsilon}\right)  \left\Vert r\right\Vert _{L^{\infty}(0,t;L^{2}%
({\mathbb{T}}^{3}))}.
\]
This last estimate along with $\left(  \text{\ref{contraction_4}}\right)  $
gives us%
\[
\left\Vert \nabla\left(  B\left(  v_{1}\right)  -B(v_{2})\right)  \right\Vert
_{L^{2}((0,t)\times{\mathbb{T}}^{3})}\leq t^{\frac{1}{2}}C_{\delta
,\varepsilon}\exp\left(  \left(  1+t\right)  C_{\delta,\varepsilon}\right)
\left\Vert \nabla v_{1}-\nabla v_{2}\right\Vert _{L^{2}((0,t)\times
{\mathbb{T}}^{3})}.
\]
$\ $We conclude that for a small $T^{\star}$ the operator has a fixed point
$u\in E_{T^{\star}}$ which verifies \eqref{approx}. As the pair $\left(
\rho,u\right)  $ solution of the above system verifies by integration of the
first equation%
\[
\int_{{\mathbb{T}}^{3}}\rho\left(  t\right)  +\varepsilon\int_{0}^{t}%
\int_{{\mathbb{T}}^{3}}\rho^{2\gamma}+\varepsilon\int_{0}^{t}\int
_{{\mathbb{T}}^{3}}\rho^{2\gamma+1}+\varepsilon\int_{0}^{t}\int_{{\mathbb{T}%
}^{3}}\rho^{3}=\int_{{\mathbb{T}}^{3}}\rho_{0}^{reg},
\]
using the second equation of $\left(  \text{\ref{approx}}\right)  $ we see
that the last relation implies that%
\[
\left\Vert \nabla u\right\Vert _{L^{2}((0,T^{\star})\times{{\mathbb{T}}^{3}}%
)}\leq\tilde{C}\left(  c,\gamma\right)  \left(  \frac{1}{\varepsilon}%
{\displaystyle\int_{\mathbb{T}^{3}}}
\rho_{0}^{reg}+\left\Vert f\right\Vert _{L_{t}^{2}L^{\frac{6}{5}}}^{2}\right)
^{\frac{1}{2}}.
\]
with the same $\tilde{C}\left(  c,\gamma\right)  $ appearing in $\left(
\text{\ref{operator_Boundness}}\right)  $. Thus, we may re-iterate the fixed
point argument. This implies that the solution $\left(  \rho,u\right)  $ of
$\left(  \text{\ref{approx}}\right)  $ is global.

\subsection{The limit $\delta\rightarrow0$}

We consider $\left(  \rho^{\delta},u^{\delta}\right)  $ a sequence of
solutions to%

\begin{equation}
\left\{
\begin{array}
[c]{l}%
\partial_{t}\rho^{\delta}+\operatorname{div}\left(  \rho^{\delta}%
\omega_{\delta}\ast u^{\delta}\right)  =\varepsilon\Delta\rho^{\delta
}-\varepsilon(\rho^{\delta})^{2\gamma}-\varepsilon(\rho^{\delta})^{2\gamma
+1}-\varepsilon\left(  \rho^{\delta}\right)  ^{3},\\
\mathcal{A}u^{\delta}+\nabla\omega_{\delta}\ast\left(  \rho^{\delta}\right)
^{\gamma}=f,\\
\rho_{|t=0}=\omega_{\delta}\ast\rho_{0}%
\end{array}
\right.  \tag{$\mathcal{S}_\varepsilon,\delta$}\label{approximation_level3}%
\end{equation}
The sequence verifies the following estimates uniformly in $\delta:$%
\begin{equation}
\left\{
\begin{array}
[c]{l}%
\displaystyle\int_{{\mathbb{T}}^{3}}\rho^{\delta}\left(  t\right)
+\varepsilon\int_{0}^{t}\int_{{\mathbb{T}}^{3}}\left(  \rho^{\delta}\right)
^{2\gamma}+\varepsilon\int_{0}^{t}\int_{{\mathbb{T}}^{3}}\left(  \rho^{\delta
}\right)  ^{2\gamma+1}+\varepsilon\int_{0}^{t}\int_{{\mathbb{T}}^{3}}\left(
\rho^{\delta}\right)  ^{3}=\int_{{\mathbb{T}}^{3}}\omega_{\delta}\ast\rho
_{0}\leq\int_{\mathbb{T}^{3}}\rho_{0},\\
\displaystyle\int_{{\mathbb{T}}^{3}}\left(  \rho^{\delta}\right)  ^{\gamma
}\left(  t\right)  +\frac{c\left(  \gamma-1\right)  }{2}\int_{0}^{t}%
\int_{{\mathbb{T}}^{3}}\left\vert \nabla u^{\delta}\right\vert ^{2}\\
\displaystyle\hskip3cm+\varepsilon\gamma\int_{0}^{t}\int_{{\mathbb{T}}^{3}%
}\left(  \rho^{\delta}\right)  ^{3\gamma-1}+\varepsilon\gamma\int_{0}^{t}%
\int_{{\mathbb{T}}^{3}}\left(  \rho^{\delta}\right)  ^{3\gamma}+\varepsilon
\gamma\int_{0}^{t}\int_{{\mathbb{T}}^{3}}\left(  \rho^{\delta}\right)
^{\gamma+2}\\
\displaystyle\hskip5cm+4\varepsilon\lbrack1-\frac{1}{\gamma}]\int_{0}^{t}%
\int_{{\mathbb{T}}^{3}}\left\vert \nabla\left(  \rho^{\delta}\right)
^{\frac{\gamma}{2}}\right\vert \leq C\left(  c,\gamma\right)  \left(
\int_{{\mathbb{T}}^{3}}\rho_{0}^{\gamma}+\left\Vert f\right\Vert _{L_{t}%
^{2}L^{\frac{6}{5}}}^{2}\right)  ,\\
\displaystyle\left\Vert \omega_{\delta}\ast\left(  \rho^{\delta}\right)
^{\gamma}\right\Vert _{L^{2}((0,T)\times{\mathbb{T}}^{3})}\leq\sqrt{t}%
{\displaystyle\int_{\mathbb{T}^{3}}}
\rho^{\gamma}+\left\Vert \Delta^{-1}\operatorname{div}\mathcal{A}u^{\delta
}\right\Vert _{L^{2}((0,T)\times{\mathbb{T}}^{3})}+\left\Vert \Delta
^{-1}\operatorname{div}f\right\Vert _{L^{2}((0,T)\times{\mathbb{T}}^{3})}\\
\leq C\left(  \gamma,c\right)  \left(  \sqrt{t}+\max\left\{  1,\left\Vert
A\right\Vert _{L^{\infty}}\right\}  \right)  \left(  \int_{{\mathbb{T}}^{3}%
}\rho_{0}^{\gamma}+\left\Vert f\right\Vert _{L_{t}^{2}L^{\frac{6}{5}}}%
^{2}\right)  ^{\frac{1}{2}}.
\end{array}
\right.  \label{uniform_in_delta}%
\end{equation}
Moreover, we have that%
\begin{align*}
&  \frac{1}{2}\int_{{\mathbb{T}}^{3}}\left(  \rho^{\delta}\right)
^{2}+\varepsilon\int_{0}^{t}\int_{{\mathbb{T}}^{3}}\left\vert \nabla
\rho^{\delta}\right\vert ^{2}\\
&  +\varepsilon\int_{0}^{t}\int_{{\mathbb{T}}^{3}}\left(  \rho^{\delta
}\right)  ^{2\gamma+1}+\varepsilon\int_{0}^{t}\int_{{\mathbb{T}}^{3}}\left(
\rho^{\delta}\right)  ^{2\gamma+2}+\varepsilon\int_{0}^{t}\int_{{\mathbb{T}%
}^{3}}\left(  \rho^{\delta}\right)  ^{4}=\gamma\int_{0}^{t}\int_{{\mathbb{T}%
}^{3}}\left(  \rho^{\delta}\right)  ^{2}\operatorname{div}\left(
\omega_{\delta}\ast u^{\delta}\right) \\
&  \hskip6cm\leq\frac{\varepsilon}{2}\int_{0}^{t}\int_{{\mathbb{T}}^{3}%
}\left(  \rho^{\delta}\right)  ^{4}+\frac{\gamma^{2}}{2\varepsilon}\int
_{0}^{t}\int_{{\mathbb{T}}^{3}}\left(  \omega_{\delta}\ast\operatorname{div}%
u^{\delta}\right)  ^{2}%
\end{align*}
and owing to the uniform bound on $\nabla u^{\delta}$ ensured by the estimates
$\left(  \text{\ref{uniform_in_delta}}\right)  $ we get that%
\begin{align}
&  \frac{1}{2}\int_{{\mathbb{T}}^{3}}\left(  \rho^{\delta}\right)
^{2}+\varepsilon\int_{0}^{t}\int_{{\mathbb{T}}^{3}}\left\vert \nabla
\rho^{\delta}\right\vert ^{2}+\varepsilon\int_{0}^{t}\int_{{\mathbb{T}}^{3}%
}\left(  \rho^{\delta}\right)  ^{2\gamma+1}+\varepsilon\int_{0}^{t}%
\int_{{\mathbb{T}}^{3}}\left(  \rho^{\delta}\right)  ^{2\gamma+2}%
+\frac{\varepsilon}{2}\int_{0}^{t}\int_{{\mathbb{T}}^{3}}\left(  \rho^{\delta
}\right)  ^{4}\nonumber\\
&  \leq\frac{C\left(  \gamma,\left\Vert A\right\Vert _{L^{\infty}}\right)
}{\varepsilon}\left(  \int_{{\mathbb{T}}^{3}}\rho_{0}^{\gamma}+\left\Vert
f\right\Vert _{L_{t}^{2}L^{\frac{6}{5}}}^{2}\right)  .
\label{uniform_in_delta2}%
\end{align}
Moreover, we have that
\begin{equation}
\partial_{t}\rho^{\delta}\text{ is bounded uniformly in }W^{-1,1}\left(
\left(  0,T\right)  \times L^{1}\left(  \mathbb{T}^{3}\right)  \right)
+L^{1}\left(  \left(  0,T\right)  \times\mathbb{T}^{3}\right)
\label{partial_t_uniform_in_delta}%
\end{equation}
The estimates $\left(  \text{\ref{uniform_in_delta}}\right)  $, $\left(
\text{\ref{uniform_in_delta2}}\right)  $ and $\left(
\text{\ref{partial_t_uniform_in_delta}}\right)  $ are enough in order to pass
to the limit when $\delta\rightarrow0$ such that we obtain the existence of a
solution of system
\[
\left\{
\begin{array}
[c]{l}%
\partial_{t}\rho+\operatorname{div}\left(  \rho u\right)  =\varepsilon
\Delta\rho-\varepsilon\rho^{2\gamma}-\varepsilon\rho^{2\gamma+1}%
-\varepsilon\rho^{3},\\
\mathcal{A}u+\nabla\rho^{\gamma}=f,\\
\rho_{|t=0}=\rho_{0}%
\end{array}
\right.
\]
which verifies the following bounds%
\begin{equation}
\left\{
\begin{array}
[c]{l}%
\displaystyle\int_{{\mathbb{T}}^{3}}\rho\left(  t\right)  +\varepsilon\int
_{0}^{t}\int_{{\mathbb{T}}^{3}}\rho^{2\gamma}+\varepsilon\int_{0}^{t}%
\int_{{\mathbb{T}}^{3}}\rho^{2\gamma+1}+\varepsilon\int_{0}^{t}\int
_{{\mathbb{T}}^{3}}\rho^{3}=\int_{{\mathbb{T}}^{3}}\rho_{0},\\
\displaystyle\int_{{\mathbb{T}}^{3}}\rho\left(  t\right)  +\frac{c\left(
\gamma-1\right)  }{2}\int_{0}^{t}\int_{{\mathbb{T}}^{3}}\left\vert
u\right\vert ^{2}\\
\displaystyle\hskip3cm+\varepsilon\gamma\int_{0}^{t}\int_{{\mathbb{T}}^{3}%
}\rho^{3\gamma-1}+\varepsilon\gamma\int_{0}^{t}\int_{{\mathbb{T}}^{3}}%
\rho^{3\gamma}+\varepsilon\gamma\int_{0}^{t}\int_{{\mathbb{T}}^{3}}%
\rho^{\gamma+2}\\
\displaystyle\hskip5cm+4\varepsilon\lbrack1-\frac{1}{\gamma}]\int_{0}^{t}%
\int_{{\mathbb{T}}^{3}}\left\vert \nabla\rho^{\frac{\gamma}{2}}\right\vert
\leq C\left(  c,\gamma\right)  \left(  \int_{{\mathbb{T}}^{3}}\rho_{0}%
^{\gamma}+\left\Vert f\right\Vert _{L_{t}^{2}L^{\frac{6}{5}}}^{2}\right)  ,\\
\displaystyle\left\Vert \rho^{\gamma}\right\Vert _{L^{2}((0,T)\times
{\mathbb{T}}^{3})}\leq C\left(  \gamma,c\right)  \left(  \sqrt{t}+\max\left\{
1,\left\Vert A\right\Vert _{L^{\infty}}\right\}  \right)  \left(
\int_{{\mathbb{T}}^{3}}\rho_{0}^{\gamma}+\left\Vert f\right\Vert _{L_{t}%
^{2}L^{\frac{6}{5}}}^{2}\right)  ^{\frac{1}{2}}.
\end{array}
\right.
\end{equation}

\subsection{Weak stability result for the perturbed system with diffusion and
drag terms}

In view of what was proved in the last section, let us consider a sequence
$\left(  \rho^{\varepsilon},u^{\varepsilon}\right)  $ of solutions of
\begin{equation}
\left\{
\begin{array}
[c]{l}%
\partial_{t}\rho^{\varepsilon}+\operatorname{div}\left(  \rho^{\varepsilon
}u^{\varepsilon}\right)  =\varepsilon\Delta\rho^{\varepsilon}-\varepsilon
\left(  \rho^{\varepsilon}\right)  ^{2\gamma}-\varepsilon(\rho^{\varepsilon
})^{3},\\
\mathcal{A}u^{\varepsilon}+\nabla\left(  \rho^{\varepsilon}\right)  ^{\gamma
}=f,\\
\rho_{|t=0}=\rho_{0}%
\end{array}
\right.  \tag{$\mathcal{S}_\varepsilon$}\label{system_eps}%
\end{equation}
which verifies the following estimates uniformly in $\varepsilon$%
\begin{equation}
\left\{
\begin{array}
[c]{l}%
\displaystyle\int_{{\mathbb{T}}^{3}}\rho^{\varepsilon}\left(  t\right)
+\varepsilon\int_{0}^{t}\int_{{\mathbb{T}}^{3}}\left(  \rho^{\varepsilon
}\right)  ^{2\gamma}+\varepsilon\int_{0}^{t}\int_{{\mathbb{T}}^{3}}\left(
\rho^{\varepsilon}\right)  ^{2\gamma+1}+\varepsilon\int_{0}^{t}\int
_{{\mathbb{T}}^{3}}\left(  \rho^{\varepsilon}\right)  ^{3}=\int_{{\mathbb{T}%
}^{3}}\rho_{0},\\
\displaystyle\int_{{\mathbb{T}}^{3}}\left(  \rho^{\varepsilon}\right)
^{\gamma}\left(  t\right)  +\left(  \gamma-1\right)  \int_{0}^{t}%
\int_{{\mathbb{T}}^{3}}\tau^{\varepsilon}:\nabla u^{\varepsilon}\\
\displaystyle\hskip2cm+\varepsilon\gamma\int_{0}^{t}\int_{{\mathbb{T}}^{3}%
}\left(  \rho^{\varepsilon}\right)  ^{3\gamma-1}+\varepsilon\gamma\int_{0}%
^{t}\int_{{\mathbb{T}}^{3}}\left(  \rho^{\varepsilon}\right)  ^{3\gamma
}+\varepsilon\gamma\int_{0}^{t}\int_{{\mathbb{T}}^{3}}\left(  \rho
^{\varepsilon}\right)  ^{\gamma+2}\\
\displaystyle\hskip4cm+4\varepsilon\lbrack1-\frac{1}{\gamma}]\int_{0}^{t}%
\int_{{\mathbb{T}}^{3}}\left\vert \nabla\left(  \rho^{\varepsilon}\right)
^{\frac{\gamma}{2}}\right\vert ^{2}\leq C\left(  \gamma,c\right)  \left(
\int_{{\mathbb{T}}^{3}}\rho_{0}^{\gamma}+\left\Vert f\right\Vert _{L_{t}%
^{2}L^{\frac{6}{5}}}^{2}\right)  ,\\
\displaystyle\left\Vert \left(  \rho^{\varepsilon}\right)  ^{\gamma
}\right\Vert _{L^{2}((0,T)\times{\mathbb{T}}^{3})}\leq C\left(  \gamma
,c\right)  \left(  \sqrt{t}+\max\left\{  1,\left\Vert A\right\Vert
_{L^{\infty}}\right\}  \right)  \left(  \int_{{\mathbb{T}}^{3}}\rho
_{0}^{\gamma}+\left\Vert f\right\Vert _{L_{t}^{2}L^{\frac{6}{5}}}^{2}\right)
^{\frac{1}{2}}.
\end{array}
\right.  \label{uniform_in_eps}%
\end{equation}
In the following we show that it is possible to slightly modify the proof of
stability in order to show that the limiting function $\left(  \rho,u\right)
$ is a solution of the semi-stationary Stokes system. Indeed, let us observe
that%
\begin{align*}
&  \gamma\left(  h+\omega_{\varepsilon^{\prime}}\ast\left(  \rho^{\varepsilon
}\right)  \right)  ^{\gamma-1}\Delta\omega_{\varepsilon^{\prime}}\ast\left(
\rho^{\varepsilon}\right) \\
&  =\Delta\left(  \left(  h+\omega_{\varepsilon^{\prime}}\ast\left(
\rho^{\varepsilon}\right)  \right)  ^{\gamma}\right)  -\nabla\left(
h+\omega_{\varepsilon^{\prime}}\ast\left(  \rho^{\varepsilon}\right)  \right)
^{\gamma-1}\nabla\omega_{\varepsilon^{\prime}}\ast\left(  \rho^{\varepsilon
}\right) \\
&  =\Delta\left(  \left(  h+\omega_{\varepsilon^{\prime}}\ast\left(
\rho^{\varepsilon}\right)  \right)  ^{\gamma}\right)  -\left(  \gamma
-1\right)  \left(  h+\omega_{\varepsilon^{\prime}}\ast\left(  \rho
^{\varepsilon}\right)  \right)  ^{\gamma-2}\nabla\omega_{\varepsilon^{\prime}%
}\ast\left(  \rho^{\varepsilon}\right)  \nabla\omega_{\varepsilon^{\prime}%
}\ast\left(  \rho^{\varepsilon}\right) \\
&  =\Delta\left(  \left(  h+\omega_{\varepsilon^{\prime}}\ast\left(
\rho^{\varepsilon}\right)  \right)  ^{\gamma}\right)  -\gamma\frac{\left(
\gamma-1\right)  }{\left(  \frac{\gamma}{2}\right)  ^{2}}\nabla\left(
h+\omega_{\varepsilon^{\prime}}\ast\left(  \rho^{\varepsilon}\right)  \right)
^{\frac{\gamma}{2}}\nabla\left(  h+\omega_{\varepsilon^{\prime}}\ast\left(
\rho^{\varepsilon}\right)  \right)  ^{\frac{\gamma}{2}}.
\end{align*}
Thus, in the sense of distributions, we get that%
\[
\gamma\left(  h+\omega_{\varepsilon^{\prime}}\ast\left(  \rho^{\varepsilon
}\right)  \right)  ^{\gamma-1}\Delta\omega_{\varepsilon^{\prime}}\ast\left(
\rho^{\varepsilon}\right)  \underset{\varepsilon^{\prime},h\rightarrow
0}{\rightarrow}\Delta\left(  \rho^{\varepsilon}\right)  ^{\gamma}%
-4\,[1-\frac{1}{\gamma}]\,\left\vert \nabla\left(  \rho^{\varepsilon}\right)
^{\frac{\gamma}{2}}\right\vert ^{2}.
\]
Also, we have that%
\[
\left\{
\begin{array}
[c]{c}%
\left(  h+\omega_{\varepsilon^{\prime}}\ast\left(  \rho^{\varepsilon}\right)
\right)  ^{\gamma-1}\omega_{\varepsilon^{\prime}}\ast\left(  \rho
^{\varepsilon}\right)  ^{2\gamma}\underset{\varepsilon^{\prime},h\rightarrow
0}{\rightarrow}(\rho^{\varepsilon})^{3\gamma-1}\text{ in }L_{t,x}^{1},\\
\left(  h+\omega_{\varepsilon^{\prime}}\ast\left(  \rho^{\varepsilon}\right)
\right)  ^{\gamma-1}\omega_{\varepsilon^{\prime}}\ast\left(  \rho
^{\varepsilon}\right)  ^{2\gamma+1}\underset{\varepsilon^{\prime}%
,h\rightarrow0}{\rightarrow}(\rho^{\varepsilon})^{3\gamma}\text{ in }%
L_{t,x}^{1}\\
\left(  h+\omega_{\varepsilon^{\prime}}\ast\left(  \rho^{\varepsilon}\right)
\right)  ^{\gamma-1}\omega_{\varepsilon^{\prime}}\ast\left(  \rho
^{\varepsilon}\right)  ^{3}\underset{\varepsilon^{\prime},h\rightarrow
0}{\rightarrow}(\rho^{\varepsilon})^{\gamma+2}\text{ in }L_{t,x}^{1}.
\end{array}
\right.
\]
We may thus write the renormalized equation for $\left(  \rho^{\varepsilon
}\right)  ^{\gamma}$ in two ways. First, we have that%
\begin{align*}
&  \partial_{t}(\rho^{\varepsilon})^{\gamma}+\operatorname{div}\left(
(\rho^{\varepsilon})^{\gamma}u^{\varepsilon}\right)  +\left(  \gamma-1\right)
(\rho^{\varepsilon})^{\gamma}\operatorname{div}u^{\varepsilon}\\
&  =\varepsilon\Delta\left(  \rho^{\varepsilon}\right)  ^{\gamma}%
-4\varepsilon\,[1-\frac{1}{\gamma}]\,\left\vert \nabla\left(  \rho
^{\varepsilon}\right)  ^{\frac{\gamma}{2}}\right\vert ^{2}-\varepsilon
\,(\rho^{\varepsilon})^{3\gamma-1}-\varepsilon\,(\rho^{\varepsilon})^{3\gamma
}-\varepsilon\,(\rho^{\varepsilon})^{\gamma+2}.
\end{align*}
which we will use to obtain uniform bounds for $\left(  \partial
_{t}u^{\varepsilon}\right)  _{\varepsilon>0}$. Secondly, we have that%
\begin{align*}
&  \partial_{t}(\rho^{\varepsilon})^{\gamma}+\gamma\operatorname{div}\left(
(\rho^{\varepsilon})^{\gamma}u^{\varepsilon}\right) \\
&  =\left(  \gamma-1\right)  \operatorname{div}(u^{\varepsilon}\tau
^{\varepsilon})-\left(  \gamma-1\right)  \tau^{\varepsilon}:\nabla
u^{\varepsilon}+u^{\varepsilon}f\\
&  +\varepsilon\Delta\left(  \rho^{\varepsilon}\right)  ^{\gamma}%
-4\varepsilon\,[1-\frac{1}{\gamma}]\,\left\vert \nabla\left(  \rho
^{\varepsilon}\right)  ^{\frac{\gamma}{2}}\right\vert ^{2}-\varepsilon
\,(\rho^{\varepsilon})^{3\gamma-1}-\varepsilon\,(\rho^{\varepsilon})^{3\gamma
}-\varepsilon\,(\rho^{\varepsilon})^{\gamma+2}.
\end{align*}
which is used for the compactness argument.

Let us observe that the time derivative of $u$ verifies%
\begin{align*}
{\mathcal{A}}\partial_{t}u^{\varepsilon}  &  =\mathrm{div}(\partial
_{t}A(t,x)D(u^{\varepsilon}))+\partial_{t}f-\nabla\partial_{t}(\rho
^{\varepsilon})^{\gamma}\\
&  =\mathrm{div}(\partial_{t}A(t,x)D(u^{\varepsilon}))+\partial_{t}f\\
&  \hskip.5cm-\nabla\operatorname{div}(\left(  \rho^{\varepsilon})^{\gamma
}u^{\varepsilon}\right)  -\left(  \gamma-1\right)  \nabla\left(
(\rho^{\varepsilon})^{\gamma}\operatorname{div}u^{\varepsilon}\right)
-\varepsilon\nabla\Delta\left(  \rho^{\varepsilon}\right)  ^{\gamma}\\
&  \hskip.5cm+4\varepsilon\,[1-\frac{1}{\gamma}]\nabla\,\left\vert
\nabla\left(  \rho^{\varepsilon}\right)  ^{\frac{\gamma}{2}}\right\vert
^{2}+\varepsilon\nabla\,(\rho^{\varepsilon})^{3\gamma-1}+\varepsilon
\nabla\,(\rho^{\varepsilon})^{3\gamma}+\varepsilon\nabla\,(\rho^{\varepsilon
})^{\gamma+2}.
\end{align*}
Also, we have that%
\[
\varepsilon\nabla\left(  \rho^{\varepsilon}\right)  ^{\gamma}=2\varepsilon
\left(  \rho^{\varepsilon}\right)  ^{\frac{\gamma}{2}}\nabla\left(
\rho^{\varepsilon}\right)  ^{\frac{\gamma}{2}},
\]
such that we obtain%
\begin{align*}
\varepsilon\int_{0}^{t}\left\Vert \nabla\left(  \rho^{\varepsilon}\right)
^{\gamma}\right\Vert _{L^{\frac{3}{2}}}  &  \leq\varepsilon\int_{0}^{t}\left(
\int_{\mathbb{T}^{3}}\left(  \rho^{\varepsilon}\right)  ^{3\gamma}\right)
^{\frac{1}{3}}\left\Vert \nabla\left(  \rho^{\varepsilon}\right)
^{\frac{\gamma}{2}}\right\Vert _{L^{2}}\\
&  \leq C\left(  t^{\frac{1}{6}}+\varepsilon\int_{0}^{t}\int_{\mathbb{T}^{3}%
}\left(  \rho^{\varepsilon}\right)  ^{3\gamma}+\varepsilon\int_{0}^{t}%
\int_{\mathbb{T}^{3}}\left\vert \nabla\left(  \rho^{\varepsilon}\right)
^{\frac{\gamma}{2}}\right\vert ^{2}\right)
\end{align*}
and we see that $\left(  \nabla\left(  \rho^{\varepsilon}\right)  ^{\gamma
}\right)  _{\varepsilon>0}$ is uniformly bounded in $L_{t}^{1}(L^{\frac{3}{2}%
}\left(  \mathbb{T}^{3}\right)  )$. It remains to write that
\[
{\mathcal{A}}\partial_{t}u^{\varepsilon}={\mathcal{A}}\phi_{1}^{\varepsilon
}+{\mathcal{A}}\phi_{2}^{\varepsilon}+{\mathcal{A}}\phi_{3}^{\varepsilon},
\]
with
\[
\left\{
\begin{array}
[c]{l}%
{\mathcal{A}}\phi_{1}^{\varepsilon}=\mathrm{div}(\partial_{t}%
A(t,x)D(u^{\varepsilon})),\\
{\mathcal{A}}\phi_{2}^{\varepsilon}=-\nabla\left\{  \operatorname{div}(\left(
\rho^{\varepsilon})^{\gamma}u^{\varepsilon}\right)  +\left(  \gamma-1\right)
(\rho^{\varepsilon})^{\gamma}\operatorname{div}u^{\varepsilon}+\varepsilon
\Delta\left(  \rho^{\varepsilon}\right)  ^{\gamma}\right\}  ,\\
{\mathcal{A}}\phi_{3}^{\varepsilon}=\nabla\left\{  4\varepsilon\,[1-\frac
{1}{\gamma}]\left\vert \nabla\left(  \rho^{\varepsilon}\right)  ^{\frac
{\gamma}{2}}\right\vert ^{2}+\varepsilon(\rho^{\varepsilon})^{3\gamma
-1}+\varepsilon\nabla\,(\rho^{\varepsilon})^{3\gamma}+\varepsilon
(\rho^{\varepsilon})^{\gamma+2}\right\}  .
\end{array}
\right.
\]
Proceeding as in Proposition \ref{extra_integrability} we obtain an uniform
bound for $\left(  \partial_{t}u^{\varepsilon}\right)  _{\varepsilon>0}$ in
$L_{t}^{1}(L^{\frac{3}{2}-}\left(  \mathbb{T}^{3}\right)  )$.

Taking in consideration the renormalized equation for $\rho$, we conclude that%
\begin{align}
&  \partial_{t}\left(  \overline{\rho^{\gamma}}-\rho^{\gamma}\right)
+\operatorname{div}\left(  \left(  \overline{\rho^{\gamma}}-\rho^{\gamma
}\right)  u\right)  +\left(  \gamma-1\right)  \left(  \overline{\rho^{\gamma}%
}-\rho^{\gamma}\right)  \operatorname{div}u\nonumber\\
&  =-\left(  \gamma-1\right)  \left\{  \overline{\tau:\nabla u}-\tau:\nabla
u\right\}  -\nu\label{difference}%
\end{align}
where $\nu$ is a positive measure i.e.%
\[
\nu=\lim_{\varepsilon\rightarrow0}\Bigl(4\varepsilon\lbrack1-\frac{1}{\gamma
}]\left\vert \nabla\left(  \rho^{\varepsilon}\right)  ^{\frac{\gamma}{2}%
}\right\vert ^{2}+\varepsilon(\rho^{\varepsilon})^{3\gamma-1}+\varepsilon
(\rho^{\varepsilon})^{3\gamma}+(\rho^{\varepsilon})^{\gamma+2}\Bigr)
\]
Arguing along the same lines as in Subsection \eqref{weakstab} we obtain that
$\overline{\rho^{\gamma}}=\rho^{\gamma}$. This concludes the proof of the
existence part of Theorem \ref{Main}.

\section{Applications to other systems\label{Applications}}

The objective of this paper is to give a proof \`{a} la Lions for the problem
of existence of weak solutions for the Quasi-Stationary Stokes system. In the
presentation, we choose to keep the model as simple as possible in order to
avoid technical difficulties that would hinder the main idea to obtain
compactness for the density: comparing the limit of the energy associated to a
sequence of weak-solutions with the energy associated to the system verified
by the limit. The objective of this section is to briefly discuss some further
extensions of our work that require only slight modifications of the arguments
presented above in order to be formally proved. First of all our results apply
to any perturbation of system $\left(  \text{\ref{ANISYS}}\right)  $ in the
form:
\begin{equation}
\left\{
\begin{array}
[c]{l}%
\partial_{t}\rho+\operatorname{div}\left(  \rho u\right)  =0,\\
-\mathrm{div}\ {\tau}+a\nabla\rho^{\gamma}+Lu=f,
\end{array}
\right.  \label{modif}%
\end{equation}
where $L:\left[  L^{2}\left(  \mathbb{T}^{3}\right)  \right]  ^{3}\rightarrow$
$\left[  H^{-1}\left(  \mathbb{T}^{3}\right)  \right]  ^{3}$ is a linear
bounded operator such that
\[
\int_{\mathbb{T}^{3}}\left\langle Lu,u\right\rangle \geq0,\text{ }\partial
_{t}\left(  Lu\right)  =L\partial_{t}u
\]
for simplicity. An interesting choice that fits in this framework is
\[
\left(  Lu\right)  ^{i}=\partial_{j}\left(  \mu\ast(Du)_{ij}-\lambda
\ast\operatorname{div}u\delta_{ij}\right)
\]
where $\mu,\lambda$ are some smooth convolution kernels which amounts in
changing the stress tensor into%
\[
\tau_{ij}=\tau_{ij}^{loc}+\tau_{ij}^{nonloc}=A_{ijkl}\left[  D(u)\right]
_{kl}+\mu\ast(Du)_{ij}-\lambda\ast\operatorname{div}u\delta_{ij}%
\]
Of course, one has to assume appropriate conditions such as to ensure
coercivity. Then existence of weak-solutions for system $\left(
\text{\ref{modif}}\right)  $ follows without any significant modifications.
Nonlocal effects are important in micro-fluidics where one is interested in
fluids flowing within thin domains see for instance \cite{Erin}. Another
common choice for the operator $L$, , see \cite{Li}, modeling the effect of an
electromagnetic field on the fluid is
\[
Lu=B\times(B\times u),
\]
where $B\in L^{\infty}\left(  \mathbb{T}^{3}\right)  $ with $B$ non-constant,
case in which we can incorporate also a force term of the type $\rho g$.

Another situation where the weak-stability part of our result can be adapted
without too much of an effort is given by the following stationary system
\[
\left\{
\begin{array}
[c]{l}%
\alpha\rho+\operatorname{div}\left(  \rho u\right)  =f,\\
\beta\rho u+\operatorname{div}\left(  \rho u\otimes u\right)
-\operatorname{div}\tau+a\nabla\rho^{\gamma}=g,
\end{array}
\right.
\]
where $a,\alpha,\beta>0$, $f\geq0$ and $\tau$ is as above. This later system
can be viewed as an implicit time discretization of the Navier-Stokes system.
Obviously, on may add nonlocality into the model. Note however that our
results do not apply to the case $\alpha=\beta=0$ corresponding to the
stationary Navier-Stokes system. This will be the object of a forthcoming
paper \cite{BrBu2}.

\section*{Appendix : Fourier Analysis on the torus and elliptic estimates}

In the following lines we present some results from Fourier analysis in the
periodic setting. The proofs are essentially the same as those in the whole
space presented in the book by H. Bahouri, J.-Y. Chemin, R. Danchin
\cite{BCD}, Chapter $2$ pages $52$-$53$. To simplify the presentation, assume
that $u\in L^{1}\left(  \mathbb{T}^{d}\right)  $. We start by reminding the
definition and properties of Fourier coefficients of $u$:
\[
\hat{u}_{\eta}=%
{\displaystyle\int\limits_{\mathbb{T}^{n}}}
\exp\left(  -2\pi y\cdot\eta\right)  u\left(  y\right)  dy.
\]
We recall the existence of two positive functions $\left(  \chi,\phi\right)
\in\mathcal{D}\left(  \mathbb{R}^{d}\right)  $ such that $\operatorname{Supp}%
\chi\subset B\left(  0,\frac{2}{3}\right)  $, $\operatorname{Supp}\phi
\subset\left\{  x:\frac{3}{4}\leq\left\vert x\right\vert \leq\frac{8}%
{3}\right\}  $ with the property that%
\[
\chi\left(  \eta\right)  +\sum_{j\geq-1}\phi\left(  2^{-j}\eta\right)
=1\text{ }\forall\eta\in\mathbb{T}^{d}.
\]
Next, for any $u\in L^{1}\left(  \mathbb{T}^{d}\right)  $, we introduce the
$j^{th}$-dyadic block operator defined as%
\[
\Delta_{j}^{\mathrm{per}}u\left(  x\right)  =\sum_{\eta\in\mathbb{Z}^{d}}%
\phi\left(  2^{-j}\eta\right)  \hat{u}_{\eta}\exp\left(  2\pi x\cdot
\eta\right)  .
\]
This operator localizes $u$ near its frequencies of magnitude $2^{j}$. \ Using
the Poisson summation formula we see that
\[
\Delta_{j}^{\mathrm{per}}u\left(  x\right)  =\int_{\mathbb{R}^{d}}%
2^{jd}h\left(  2^{j}(x-y)\right)  u\left(  y\right)  dy
\]
where $h$ is the Fourier inverse of $\phi$.\ This last identity is useful to
show that $\Delta_{j}^{\mathrm{per}}$ maps all $L^{p}\left(  \mathbb{T}%
^{d}\right)  $ into $L^{p}\left(  \mathbb{T}^{d}\right)  $ with norm
independent of $j$ and $p$. For all $u\in L^{1}\left(  \mathbb{T}^{d}\right)
$ we have that%
\[
u=\int_{\mathbb{T}^{d}}u+\sum_{j\geq-1}\Delta_{j}^{\mathrm{per}}u
\]
at least in the sense of distributions. We infer that for any $u\in L^{p}$
with $\int_{\mathbb{T}^{d}}u=0$ we have that%
\begin{equation}
\left\Vert u\right\Vert _{L^{p}}\leq\sum_{j\geq-1}\left\Vert \Delta
_{j}^{\mathrm{per}}u\right\Vert _{L^{p}}. \label{inequality_blocks}%
\end{equation}
Next, let us recall the celebrated Bernstein lemma.

\begin{lemma}
\label{Bernstein}Consider any nonnegative integer $k$, a couple $p,q\in\left[
1,\infty\right]  ^{2}$ with $p\leq q$ and a function $u\in L^{1}\left(
\mathbb{T}^{d}\right)  $. Then, there exists a constant $C$ such that the
following inequalities hold true:%
\begin{equation}
\sup_{\left\vert \alpha\right\vert =k}\left\Vert \partial^{\alpha}\Delta
_{j}^{\mathrm{per}}u\right\Vert _{L^{q}}\leq C^{k+1}2^{jk+j\left(  \frac{d}%
{p}-\frac{d}{q}\right)  }\left\Vert \Delta_{j}^{\mathrm{per}}u\right\Vert
_{L^{p}}, \label{Bernstein1}%
\end{equation}
and
\begin{equation}
C^{-k-1}2^{jk}\left\Vert \Delta_{j}^{\mathrm{per}}u\right\Vert _{L^{p}}%
\leq\sup_{\left\vert \alpha\right\vert =k}\left\Vert \partial^{\alpha}%
\Delta_{j}^{\mathrm{per}}u\right\Vert _{L^{p}}\leq C^{k+1}2^{jk}\left\Vert
\Delta_{j}^{\mathrm{per}}u\right\Vert _{L^{p}}. \label{Bernstein2}%
\end{equation}

\end{lemma}

The following proposition will be very useful in establishing estimates for
the Poisson problem.

\begin{proposition}
\label{multiplier}Consider $m\in\mathbb{R}$ and a smooth function
$\sigma:\mathbb{R}^{d}\backslash\{0\}\rightarrow\mathbb{R}$ such that for all
multi-index $\alpha$ with $\left\vert \alpha\right\vert \leq2+2\left[
d/2\right]  $, there exists a constant $C_{\alpha}$ such that:%
\[
\forall\xi\in\mathbb{R}^{d}\backslash\{0\}\text{ }:\text{\ }\left\vert
\partial^{\alpha}\sigma\left(  \xi\right)  \right\vert \leq C_{\alpha
}\left\vert \xi\right\vert ^{m-\left\vert \alpha\right\vert }.
\]
Then for any $p\in\left[  1,\infty\right]  $ we have that%
\[
\left\Vert \sigma\left(  D\right)  \Delta_{j}^{\mathrm{per}}v\right\Vert
_{L^{p}}\leq2^{jm}\left\Vert \Delta_{j}^{\mathrm{per}}v\right\Vert _{L^{p}}%
\]
where
\[
\sigma\left(  D\right)  \Delta_{j}^{\mathrm{per}}v=\sum_{\eta\in\mathbb{Z}%
^{d}}\phi\left(  2^{-j}\eta\right)  \sigma\left(  \eta\right)  \hat{u}_{\eta
}\exp\left(  2\pi x\cdot\eta\right)  .
\]

\end{proposition}

Finally, we use the Littlewood-Paley apparatus in order to prove the following
$3D$ estimate for the Poisson problem.

\begin{theorem}
\label{Stampacchia_Gallouet}Consider $f\in L^{1}\left(  \mathbb{T}^{3}\right)
$ such that $\int_{\mathbb{T}^{3}}f=0$ and $\psi$ solution to the Poisson
problem%
\[
\left\{
\begin{array}
[c]{c}%
-\Delta\psi=f,\\
\int_{\mathbb{T}^{3}}\psi=0
\end{array}
\right.
\]
Then there exists a constant $C$ such that for any $p\in\lbrack1,\frac{3}{2})$
we have
\[
\left\Vert \nabla\psi\right\Vert _{L^{P}}\leq C\left\Vert f\right\Vert
_{L^{1}}.
\]

\end{theorem}

\noindent\textbf{Proof.} For any $l\in\overline{1,3}$ let observe that the
function $\sigma_{l}:\mathbb{R}^{3}\backslash\{0\}\rightarrow\mathbb{R}$
defined as%
\[
\sigma_{l}\left(  \xi\right)  =\frac{i\xi_{l}}{\left\vert \xi\right\vert ^{2}%
},
\]
verifies the hypothesis of Proposition \ref{multiplier}. Next, we see that for
any $\eta\in\mathbb{Z}^{d}\backslash\{0\}$ and any $l\in\overline{1,3}$ we
have that
\[
\widehat{\partial_{l}\psi}\left(  \eta\right)  =i\eta_{l}\hat{\psi}\left(
\eta\right)  =\frac{i\eta_{l}}{\left\vert \eta\right\vert ^{2}}\hat{f}\left(
\eta\right)  =\sigma_{l}\left(  \eta\right)  \hat{f}\left(  \eta\right)
\]
such that%
\[
\Delta_{j}^{\mathrm{per}}\left(  \partial_{l}\psi\right)  =\sigma_{l}\left(
D\right)  \Delta_{j}^{\mathrm{per}}f
\]
Let $p\in\lbrack1,\frac{3}{2})$. As, $\int\partial_{l}\psi=0$ using $\left(
\text{\ref{inequality_blocks}}\right)  $, Proposition \ref{multiplier} and
Bernstein's inequality, we infer that
\begin{align*}
\left\Vert \partial_{l}\psi\right\Vert _{L^{p}}  &  \leq\sum_{j\geq
-1}\left\Vert \Delta_{j}^{\mathrm{per}}\partial_{l}\psi\right\Vert _{L^{p}%
}=\sum_{j\geq-1}\left\Vert \sigma_{l}\left(  D\right)  \Delta_{j}%
^{\mathrm{per}}f\right\Vert _{L^{p}}\leq2^{-j}\sum_{j\geq-1}\left\Vert
\Delta_{j}^{\mathrm{per}}f\right\Vert _{L^{p}}\\
&  \leq\sum_{j\geq-1}2^{j\left(  2-\frac{3}{p}\right)  }\left\Vert \Delta
_{j}^{\mathrm{per}}f\right\Vert _{L^{1}}\leq\left\Vert f\right\Vert _{L^{1}%
}\sum_{j\geq-1}2^{j\left(  2-\frac{3}{p}\right)  },
\end{align*}
where, of course the fact that $p\in\lbrack1,3/2)$ ensures the convergence of
the series $\sum_{j\geq-1}2^{j\left(  2-\frac{3}{p}\right)  }$. With this
remark we conclude the proof of Theorem \ref{Stampacchia_Gallouet}.

\begin{remark}
In fact, a more careful analysis of the proof of Theorem
\ref{Stampacchia_Gallouet} yields the following refined estimate%
\[
\left\Vert \nabla\psi\right\Vert _{L^{p}(\mathbb{T}^{3})}\lesssim\left\Vert
\nabla\psi\right\Vert _{B_{p,1}^{0}(\mathbb{T}^{3})}\lesssim\left\Vert
f\right\Vert _{B_{1,\infty}^{0}}%
\]
which is stronger than the classical result as the space of bounded measures
is continuously included in $B_{1,\infty}^{0}$.
\end{remark}

\noindent\textbf{Acknowledgments.} D. Bresch and C. Burtea are supported by
the SingFlows project, grant ANR-18-CE40-0027 and D. Bresch is also supported
by the Fraise project, grant ANR-16-CE06-0011 of the French National Research
Agency (ANR).

\end{document}